\documentclass[a4paper,11pt]{article}
 \def\draft{0}
 \if\draft1
 \usepackage[top=3cm, bottom=4.5cm, right=3.5cm, left=1.5cm]{geometry}
  \usepackage{refcheck}
 \else
\usepackage[top=3cm, bottom=2.5cm, right=2.5cm, left=2.5cm]{geometry}
 \fi
 
\usepackage[]{hyperref} 
\usepackage[utf8]{inputenc}
\usepackage[T1]{fontenc}     
\usepackage{lmodern}         

\usepackage{amsmath}
\usepackage{amsthm}
\usepackage{amssymb}
\usepackage{mathrsfs}
\usepackage[sfdefault=cmbr,OMLmathsans]{isomath}  
\usepackage{latexsym}
\usepackage{upgreek}  
\usepackage{tensor}  
\usepackage{bbm}

\usepackage{algorithm}
\usepackage{algorithmicx}
\usepackage{algpseudocode}
\usepackage{graphicx}
\usepackage{subfig} 
\usepackage{xifthen}
\usepackage{rotating}

\usepackage{enumerate}
\usepackage[dvipsnames]{xcolor}
\usepackage{multirow}
\usepackage{soul}
\usepackage{diagbox} 

\newcommand{\comment}[1]{}

\newcommand{\MBC}{\ensuremath{\B{\M{C}}}} 
\newcommand{\MBb}{\ensuremath{\B{\M{b}}}}
\newcommand{\tMBC}{\widetilde{\B{\M{C}}}} 
\newcommand{\tMBb}{\widetilde{\B{\M{b}}}}
\newcommand{\bas}[1]{\ensuremath{\B{\M{b}}}^{(#1)}}
\newcommand{\basc}[1]{\ensuremath{\M{b}}^{(#1)}}
\newcommand{\hbas}[1]{\ensuremath{\widehat{\B{\M{b}}}}^{(#1)}}

\newcommand{\ga}{\ensuremath{\textrm{Ga}}}
\newcommand{\gani}{\ensuremath{\textrm{GaNi}}}

\newcommand{\tucker}{\ensuremath{\textrm{Tucker}}}
\newcommand{\TT}{\ensuremath{\textrm{TT}}}

\newcommand{\squar}{\ensuremath{\square}}
\newcommand{\stoch}{\ensuremath{\mathrm{S}}}

\newcommand{\eg}{e.g.\ }

\newcommand{\FT}{\mathcal{F}}
\newcommand{\FFT}[1]{\FT_{#1}}
\newcommand{\iFFT}[1]{\FT_{#1}^{-1}}
\newcommand{\FFTd}[1]{\bs{\FT}_{#1}}
\newcommand{\iFFTd}[1]{\bs{\FT}_{#1}^{-1}}
\newcommand{\Tr}{\mathcal{T}}
\newcommand{\Vr}{\ensuremath{\V{r}}}

\newcommand{\B}[1]{\boldsymbol{#1}} 
\newcommand{\TA}{\ensuremath{\T{A}}}
\newcommand{\MA}{\ensuremath{\M{A}}}
\newcommand{\MBA}{\ensuremath{\B{\M{A}}}}
\newcommand{\Hnp}{H^1_{\zmean}}
\newcommand{\zmean}{0}
\newcommand{\Vw}{\ensuremath{\V{w}}}
\newcommand{\VE}{\ensuremath{\V{E}}}
\newcommand{\eff}{\mathrm{H}}
\newcommand{\VN}{\ensuremath{{\V{N}}}}
\newcommand{\tVN}{\ensuremath{{2\V{N}-\V{1}}}}
\newcommand{\pVN}{|\VN|_{\Pi}}
\newcommand{\Ve}{\ensuremath{\V{e}}}

\newcommand{\MBF}{\ensuremath{\B{\M{F}}}}

\newcommand{\MBe}{\ensuremath{\B{\M{e}}}}
\newcommand{\MBE}{\ensuremath{\B{\M{E}}}}

\newcommand{\Mv}{\ensuremath{\M{v}}}

\newcommand{\MBu}{\ensuremath{\B{\M{u}}}}
\newcommand{\MBr}{\ensuremath{\B{\M{r}}}}
\newcommand{\MBv}{\ensuremath{\B{\M{v}}}}

\newcommand{\MBw}{\ensuremath{\B{\M{w}}}}

\newcommand{\hMv}{\ensuremath{\widehat{\M{v}}}}

\newcommand{\hMBu}{\ensuremath{\widehat{\B{\M{u}}}}}
\newcommand{\hMBv}{\ensuremath{\widehat{\B{\M{v}}}}}
\newcommand{\hMBw}{\ensuremath{\widehat{\B{\M{w}}}}}
\newcommand{\MBP}{\ensuremath{\B{\M{P}}}}
\newcommand{\hMBP}{\ensuremath{\widehat{\B{\M{P}}}}}

\newcommand{\imu}{\mathrm{i}}		

\newcommand{\ZNd}{\set{Z}_{\VN}}
\newcommand{\ZtNd}{\set{Z}_{\tVN}}

\newcommand{\cT}{\mathscr{T}}

\newcommand{\cTN}{\cT_\VN}

\newcommand{\sV}{\mathcal{V}}

\newcommand{\Aeff}{A_{\eff}}

\newcommand{\puc}{\mathcal{Y}}

\newcommand{\Vv}{{\V{v}}}

\newcommand{\Vk}{{\V{k}}}
\newcommand{\Vl}{{\V{l}}}
\newcommand{\Vm}{{\V{m}}}

\newcommand{\skipifemptyarg}[1]{\ifthenelse{\isempty{#1}}{}{\left[#1\right]}}
\newcommand{\skipifscalar}[1]{\ifthenelse{\isempty{#1}}{}{;#1}}

\newcommand{\bs}[1]{\boldsymbol{#1}}
\newcommand{\scal}[2]{\bigl(#1,#2\bigr)}
\newcommand{\set}[1]{\mathbb{#1}} 
\newcommand{\M}[1]{\mathsfit{#1}} 
\newcommand{\MB}[1]{\bs{\M{#1}}} 

\newcommand{\V}[1]{\bs{#1}} 
\newcommand{\T}[1]{\bs{#1}} 
\newcommand{\ol}[1]{\overline{#1}}

\newcommand{\iter}[1]{_{(#1)}}

\newcommand{\x}{\V{x}}

\newcommand{\alp}{\ensuremath{\alpha}}

\newcommand{\del}{\ensuremath{\delta}}

\newcommand{\sR}{\set{R}}
\newcommand{\sN}{\set{N}}
\newcommand{\sZ}{\set{Z}}
\newcommand{\sC}{\set{C}}
\newcommand{\sK}{\set{K}}

\newcommand{\C}[1]{\mathcal{#1}}
\renewcommand{\tt}[1]{{\MB{#1}}}

\newcommand{\vek}[1]{\mathchoice{\displaystyle\boldsymbol#1}
{\textstyle\boldsymbol#1}{\scriptstyle\boldsymbol#1}
{\scriptscriptstyle\boldsymbol#1}}

\newcommand{\bilf}[2]{a\ifthenelse{\isempty{#1}}{}{\bigl(#1,#2\bigr)}}
\newcommand{\bilfN}[2]{a_\VN\ifthenelse{\isempty{#1}}{}{\bigl(#1,#2\bigr)}}
\providecommand{\norm}[1]{\lVert#1\rVert}

\newcommand{\D}[1]{\,{\mathrm d}#1}

\theoremstyle{plain}

\title{FFT-based homogenisation accelerated by low-rank tensor approximations}

\usepackage{authblk}
\author[1]{Jaroslav Vondřejc}
\affil[1]{Technische Universit\"{a}t Braunschweig, Institute of Scientific Computing, Mühlenpfordstrasse~23, 38106 Braunschweig, Germany,
 \href{mailto:j.vondrejc@tu-bs.de}{j.vondrejc@tu-bs.de}
 }
\affil[2]{Czech Technical University in Prague, Faculty of Civil Engineering, Thákurova~7/2077
166~29~Prague~6, Czech Republic}
\date{\today}

\author[1]{Dishi Liu}

\author[2]{Martin Ladecký}

\author[1]{Hermann G. Matthies}

\begin{document}
 \maketitle
\begin{abstract}
Fast Fourier transform (FFT) based methods have turned out to be an effective computational approach for numerical homogenisation. In particular, Fourier-Galerkin methods are computational methods for partial differential equations that are discretised with trigonometric polynomials.
Their computational effectiveness benefits from efficient FFT based algorithms as well as a favourable condition number.
Here these kind of methods are accelerated by low-rank tensor approximation techniques for a solution field using canonical polyadic, Tucker, and tensor train formats.
This reduced order model also allows to efficiently compute suboptimal global basis functions without solving the full problem.
It significantly reduces computational and memory requirements for problems with a material coefficient field that admits a moderate rank approximation.
The advantages of this approach against those using full material tensors are demonstrated using numerical examples for the model homogenisation problem that consists of a scalar linear elliptic variational problem defined in two and three dimensional settings with continuous and discontinuous heterogeneous material coefficients.
This approach opens up the potential of an efficient reduced order modelling of large scale engineering problems with heterogeneous material.
\end{abstract}


\tableofcontents

\section{Introduction}

\emph{FFT-based methods.}
A fast Fourier transform (FFT) based method has been introduced as an efficient algorithm for numerical homogenisation in 1994 by Moulinec and Suquet \cite{Moulinec1994FFT}.
The method, that has application in multiscale problems, represents an alternative discretisation approach to the finite element method.
The effectiveness of FFT-based homogenisation relies on the facts that the system matrix is never assembled, the matrix-vector product in linear iterative solvers is provided very efficiently by FFT, and the condition number is independent of discretisation parameters.

Since the seminal paper in 1994 the methodology has been significantly developed. Originally the approach has been based on Lippmann-Schwinger equation, which is a formulation incorporating Green's function for an auxiliary homogeneous problem.
Its connection to a standard variational formulation has been discovered in \cite{VoZeMa2014FFTH} by using the fact that Green's function is a projection on compatible fields (i.e.\ gradient fields in scalar elliptic problems), see \cite{Milton2002TC}.
It has allowed to fully remove the reference conductivity tensor from the formulation, and interpreted the method from the perspective of finite elements also in nonlinear problems \cite{ZeGeVoPeGe2017,GeVoZePeGe2017large}.
Moreover, the standard primal-dual variational formulations allow to compute guaranteed bounds on effective material properties \cite{VoZeMa2014GarBounds}, which provides tighter bounds than the Hashin-Shtrikman functional.

Significant attention has been focused on developing discretisation approaches that justify the original FFT-based homogenisation algorithm.
Many efforts have been made on discretisation with trigonometric polynomials, starting with \cite{ZeVoNoMa2010AFFTH} and followed by \cite{VoZeMa2014GarBounds,Vondrejc2015FFTimproved,ZeGeVoPeGe2017,Schneider2014convergence}.
Other discretisation approaches are based on pixel-wise constant basis functions
\cite{Brisard2012FFT,Brisard2010FFT}, linear hexahedral elements \cite{Schneider2016hexa}, or finite differences \cite{Willot2015,Willot2013fourier}.
The variational formulations also allowed to derive convergence of approximate solutions to the continuous one 
\cite{VoZeMa2014FFTH,Schneider2014convergence,Brisard2012FFT}.

The various discretisation approaches have been studied along with linear and non-linear solvers \cite{Eyre1999FNS,Moulinec2017,ZeVoNoMa2010AFFTH,MiVoZe2016jcp,Brisard2010FFT,Kabel2014LargeDef,ZeGeVoPeGe2017,GeVoZePeGe2017large,Schneider2016nestor}.
Other research directions focus, for example, on multiscale methods \cite{Kochmann2017,Gokuzum2017,Dietrich2017}, highly non-linear problems in solid mechanics \cite{Bertin2018,DeGeus2016d,Boeff2015,Segurado2018}, and parameter estimation features FFT and model reduction
\cite{Garcia-Cardona2017}.

\emph{Low-rank approximations.}
The general idea of low-rank approximations is to express or compress tensors with fewer parameters, which can lead to a huge reduction in requirements for computer memory and possible significant computational speed-up.
For matrices as second order tensors, the optimal low-rank approximation in mean square sense is based on the truncated singular value decomposition (SVD). A computationally cheaper choice is Cross Approximation \cite{Goreinov1997app,Bebendorf2000} which has only linear complexity in matrix size $N$.
Low-rank formats or tensors of order larger than two include the canonical polyadic (CP), Tucker, and hierarchical schemes such as the tensor train and the quantic tensor-train form of \cite{Hackbusch2014,Kolda2009}.
Low-rank formats are not only needed to compress the data tensor as the final delivered result of high-dimensional numerical modellings, but are also preferred to approximate tensors in the numerical solution process.
In \cite{Giraldi2013} the proper generalised decomposition is adopted for the construction of low-rank tensors in CP and Tucker formats in a numerical homogenisation from high-resolution images.
It is also possible to compute the tensors directly in low-rank formats, which can be provided by a suitable solver \cite{Kressner2011,Tobler2012,Dolgov2013,Ballani2013,Matthies2012}.
The rank one tensors in low-rank approximations can be seen as suboptimal global basis functions.

However, the need to compute with tensors in low-rank formats requires one to deal with operations such as addition, element-wise multiplication, or Fourier transformation. Since the low-rank tensors are described with fewer parameters, the computational complexities are typically reduced, which may lead to significant speed-up of computations.
However, performing such operations with tensors in low-rank format, it typically happens that the representation rank of the tensors grows, which calls for their truncation, i.e.\ their approximation or reparametrisation with fewer parameters while keeping a reasonable accuracy \cite{Oseledets2010,Oseledets2011,Bigoni2016}.
This truncation of tensors may be viewed as a generalisation of the rounding of numbers, which occurs when working with floating point formats.
In general, the applications of low-rank approximations are very broad, e.g.\ for stochastic problems  with high number of random parameters  \cite{Espig2014,Matthies2012,Nouy2015,Khoromskij2011}, acceleration of solutions to PDEs \cite{Khoromskij2015,Khoromskij2017}, or model order reduction \cite{Nouy2015a}, but its application to FFT-based homogenisation is new.
However, an alternative low-rank representation has been studied recently in \cite{Kochmann2018sparse}.

\emph{Structure of the paper.}
In section~\ref{sec:homog-fourier-galerkin-method}, two state-of-the-art Fourier-Galerkin methods are described for a model homogenisation problem of a scalar elliptic equation.
In particular, the two discretisation methods based on numerical and exact integration are described along with their corresponding linear systems.
Then in section~\ref{sec:fft-based-methods-with-low-rank-approximations} the low-rank approximation techniques are summarised and their application within a Fourier-Galerkin method is discussed.
In section~\ref{sec:numerical-results}, the effectiveness of low-rank approximations is demonstrated on several numerical examples.

\emph{Notation.}
We will denote vectors and matrices by boldface letters: $\V{a} = \left(a_i \right)_{i=1,2,\ldots,d}\in \sR^d$ or 
$\B{A} = (A_{ij})_{i,j=1}^d \in \sR^{d\times d}$.
Matrix-matrix and matrix-vector multiplications are denoted as $\B{C} = \B{A}\B{B}$ and $\V{c} = \B{A}\V{b}$, which in Einstein summation notation reads $C_{ik} = A_{ij}B_{jk}$ and $b_{i} = A_{ij} b_j$ respectively.
The Euclidean inner product will be referred to as $ \V{a} \cdot \V{\bar b} = \sum_i a_i \ol{b}_i$, and the induced norm as $\| \V{a} \| = \sqrt{\V{a}\cdot\V{\bar a}}$.
Vectors, matrices, and tensors such as $\MB{x}$, $\MB{b}$, and $\MBA$ arising from discretisation will be denoted by the bold sans-serif font in order to highlight their special structure.
For $\VN=(N_1,\dotsc,N_d)\in\sN^d$, the components of a tensor $\MBA\in\sR^{\VN}=\bigotimes_{\alp=1}^d \sR^{N_\alp}$  of order $d$ will be denoted as $\MA[k_1,...,k_d]$.
The multiindex notation will be also incorporated to simplify the components of the tensors, e.g.\ $\MA[k_1,...,k_d]=\MA[\V{k}]$ for a multi-index $\V{k}=[k_1,\dotsc,k_d]$.
The space $\sR^\VN$, composed of tensors of order $d$, can be considered as a vector space, which allows to talk about its dimension as the number of basis vectors, i.e. $\dim\sR^\VN=\prod_{\alp=1}^d N_\alp$. 

The space of square integrable $\puc$-periodic functions defined on a periodic cell $\puc=(-\frac{1}{2},\frac{1}{2})^d$ is denoted as $L^2(\puc)$.
The analogous space $L^2(\puc;\sR^d)$ collects $\sR^d$-valued functions $\Vv:\puc\rightarrow\sR^d$ with components $v_i$ from $L^2(\puc)$. Finally, $\Hnp(\puc)=\{v\in L^2(\puc) \,|\, \nabla v\in L^2(\puc;\sR^d), \int_{\puc} v(\x) \D{\x} = 0\}$ denotes the Sobolev space of periodic functions with zero mean.

\section{Homogenisation by Fourier-Galerkin methods}\label{sec:homog-fourier-galerkin-method}
\subsection{Model problem} 
\label{sec:problem-description-model}
%
A model problem in homogenisation \cite{Bensoussan1978book} consists of a scalar linear elliptic variational problem defined on a unit domain $\puc=(-\frac{1}{2},\frac{1}{2})^d$ in a   spatial dimension $d$  (we consider both $d=2$ and $d=3$) with material coefficients $\B{A}:\puc\rightarrow\sR^{d\times d}$, which are required to be essentially
bounded, symmetric, and uniformly elliptic. This means that for almost all $\x\in\puc$, there are constants $0 < c_A \leq C_A < +\infty$ such that
\begin{align}
\label{eq:A}
\B{A}(\x)&=\B{A}^T(\x),
&
c_A \norm{\Vv}^2
&\leq
\B{A}(\x)\Vv\cdot\Vv
\leq C_A \norm{\Vv}^2
\quad\text{ for all }\Vv \in \sR^d.
\end{align}

The homogenisation problem is focused on the computation of effective material properties $\V{A}_\eff\in\sR^{d\times d}$. Its variational formulation is based on the minimisation of a microscopic energetic functional for constant vectors $\VE\in\sR^d$, which represents an \emph{average} of the macroscopic gradient, as
\begin{equation}
\label{eq:homog}
\V{A}_\eff \VE \cdot \VE
=
\min_{v\in\Hnp(\puc)} \bilf{\VE+\nabla v}{\VE+\nabla v},
\end{equation}
where the bilinear form $a:L^2(\puc;\sR^d)\times L^2(\puc;\sR^d)\rightarrow\sR$ is defined as
\begin{align}
\label{eq:bilinear}
\bilf{\Ve}{\Vw} :=
\int_\puc \B{A}(\x) \Ve(\x) \cdot \Vw(\x) \D{\x}.
\end{align}
 The minimisation Sobolev space $\Hnp(\puc)$ consists of zero-mean $\puc$-periodic microscopic fields $v:\sR^d\rightarrow\sR$, which have locally square integrable weak gradient and finite $L^2$-norm on $\puc$;
together with \eqref{eq:A} it satisfies the existence of a unique minimiser.
Note that the minimisation problem \eqref{eq:homog} corresponds to the scalar elliptic partial differential equation $-\nabla\cdot[\TA(\x)\nabla u(\x)] = f(\x)$ with a special right-hand side $f(\x)=-\nabla\cdot\TA(\x)\VE$ and periodic boundary conditions. 

\subsection{Fourier-Galerkin methods}
\label{sec:fourier-galerkin-method}
Alternatively, the minimisers in \eqref{eq:homog} are described by a weak formulation: find $u\in \Hnp(\puc)$ such that
\begin{align}
\label{eq:WF}
\bilf{\nabla u}{\nabla v} = -\bilf{\VE}{\nabla v}
\quad\forall\, v \in \Hnp(\puc).
\end{align}
This formulation is the starting point for a discretisation using Galerkin approximations, when the trial and test spaces are substituted with finite dimensional ones.
We choose to discretise the function space using trigonometric polynomials, which leads to a Fourier-Galerkin method.

In order to compute the effective matrix $\V{A}_\eff$ one has to solve $d$ minimisation problems or weak formulation for different $\V{E}$, which are usually taken as the canonical basis of $\sR^d$.
Here we consider exclusively $\VE = (\delta_{1,i})_{i=1}^d\in\sR^d$ (i.e.\ in 3D $\VE = [1,0,0]$); therefore, the $(1,1)$-component of the homogenised properties will be of particular interest, i.e.\
$
\V{A}_\eff \VE \cdot \VE = \V{A}_{\eff,11} =: \Aeff.
$

\subsubsection{Trigonometric polynomials}
The Fourier-Galerkin method, \cite{SaVa2000PIaPDE,VoZeMa2014FFTH,Vondrejc2015FFTimproved} is built on discretisations using the space of \emph{trigonometric polynomials}
\begin{align*}
\cTN &=
\Bigl\{\sum_{\Vk\in\ZNd}{\hMv[\Vk]\varphi^{ \Vk }} \; \mid \; \hMv[\Vk]\in\sC, \text{ and }\hMv[\Vk]=\overline{\hMv[-\Vk]}\Bigr\},
\end{align*}
where $\varphi^\Vk(\x) = \exp(2\pi\imu \, \Vk\cdot\x)$ are the well-known Fourier basis functions.
The number of discretisation points $\VN=[N,\dotsc,N]\in\sR^d$ in this work take only odd values because an even $N$ introduces Nyquist frequencies that have to be omitted to obtain a conforming approximation, see \cite{VoZeMa2014GarBounds} for details.

There are also other natural basis vectors $\varphi^{\Vk}_\VN:\puc\rightarrow\sR$, the so-called fundamental trigonometric polynomials. They are expressed as a linear combination
\begin{align*}
\varphi_\VN^\Vk(\x)
=
\frac{1}{\pVN}
\sum_{\Vm \in \ZNd}
\omega_{\VN}^{-\Vk\Vm} \varphi^{\Vm}(\x)
\text{ for }
\x \in \puc,
\end{align*}
of Fourier basis function $\varphi^{\Vm}$ with complex-valued weights
$\omega_{\VN}^{\Vm\Vk} =
\exp  \bigl(2 \pi \imu\sum_{\alp=1}^d \frac{m_\alp k_\alp}{N_\alp} \bigr)$
for $\Vm,\Vk\in \ZNd$. The weights are from the discrete Fourier transform (DFT) matrices in $\sC^{\VN\times\VN}$ with components
\begin{align*}
\FFT{\VN}[\Vm,\Vk] = \frac{1}{\pVN} \omega_{\VN}^{-\Vm\Vk},
&&
\iFFT{\VN}[\Vm,\Vk] =\omega_{\VN}^{\Vm\Vk}
\quad\text{for }{\Vm,\Vk\in\ZNd}.
\end{align*}
The coefficients of trigonometric polynomials in the two different base are connected by the discrete Fourier transform (DFT), particularly expressed as
\begin{align*}
v(\x) &= \sum_{\Vk\in\ZNd}{\hMv[\Vk]\varphi^{ \Vk }}(\x)= \sum_{\Vk\in\ZNd}{\Mv[\Vk]\varphi^{ \Vk }_\VN}(\x) \quad \text{ and } \quad \hMBv=\FFT{\VN}\MBv.
\end{align*}
Due to the Dirac-delta property $\varphi_\VN^\Vl(\x_\VN^\Vk)=\delta_{\Vk\Vl}$ of the fundamental trigonometric polynomials on a regular grid of points $\x_\VN^\Vk=\frac{k_\alp}{N_\alp}$ for $\Vk,\Vl\in\ZNd$, the coefficients of the trigonometric polynomials are equal to the function values at the grid points, i.e.\ $\Mv[\Vk] = v (\x_\VN^\Vk)$.

\emph{Differential operators} are naturally applied on trigonometric polynomials. In particular the gradient
\begin{align*}
\nabla v(\x) &= \sum_{\Vk\in\ZNd} \hat{\Mv}[\Vk] \nabla\varphi^\Vk(\x) = \sum_{\Vk\in\ZNd} 2\pi\imu\Vk \hat{\Mv}[\Vk] \varphi^\Vk(\x),
\end{align*}
corresponds to the application of the operator $\widehat{\nabla}_\VN:\sC^\VN\rightarrow\sC^{d\times\VN}$ on Fourier coefficients as
$(\widehat{\nabla}_\VN\hat{\MBv})[\alp,\Vk] = 2\pi\imu k_\alp\hat{\MBv}[\Vk]$.
The adjoint operator $\widehat{\nabla}_\VN^*:\sC^{d\times\VN}\rightarrow\sC^{\VN}$ corresponding to the divergence is then expressed as
\begin{align*}
(\widehat{\nabla}_\VN^*\hat{\MBw})[\Vk] = \sum_{\alp=1}^{d} -2\pi\imu k_\alp \hat{\MBw}[\alp,\Vk]. 
\end{align*}

Then the gradient operator can be expressed with respect to the basis with fundamental trigonometric polynomials as
\begin{align}
\label{eq:grad_trig}
\nabla v(\x) &= 
 \sum_\Vk \bigl(\iFFTd{\VN}\widehat{\nabla}_\VN\FFT{\VN}\MBv\bigr)[\Vk] \varphi_\VN^\Vk(\x)
\end{align}
where the $d$-fold discrete Fourier transform (emphasises with bold)
$\FFTd{\VN} = \sC^{d\times\VN}\rightarrow\sC^{d\times\VN}$ acts individually on each component of the vector field  $(\FFTd{\VN}\MBw)[\alp] = \FFT{\VN}\MBw[\alp]$ for $\alp =1, \cdots, d$.

The numerical treatment of the weak formulation \eqref{eq:WF} or a corresponding Galerkin approximation requires the use of numerical integration. In this manuscript we incorporate two versions: an exact integration \cite{Vondrejc2015FFTimproved} as described in sub-section \ref{sec:FFTHga}, and a numerical integration as described in sub-section \ref{sec:FFTHgani}.

\subsubsection{The Fourier-Galerkin method with numerical integration (GaNi)}
\label{sec:FFTHgani}

This numerical integration based on the rectangle (or the mid-point) rule corresponds to the original Moulinec-Suquet algorithm \cite{Moulinec1994FFT,Moulinec1998NMC}, as the resulting discrete solution vectors fully coincide.
This approach, applied to the bilinear form \eqref{eq:bilinear} on regular grids, reads
\begin{align*}
\bilf{\Ve }{\Vw } \approx \bilfN{\Ve_\VN}{\Vw_\VN} = \sum_{\Vk\in\ZNd} \B{A}(\x_\VN^\Vk)\Ve_\VN(\x_\VN^\Vk)\cdot \Vw_\VN(\x_\VN^\Vk)
=\scal{\widetilde{\MBA}\MBe}{\MBw}_{\sR^{d\times\VN}},
\end{align*} 
where $\MBe$ and $\MBw$ store the function values on the grid (\eg $\M{e}[\alp,\Vk]=e_{\VN,\alp}(\x^\Vk_\VN)$), and $\widetilde{\MBA}\in\sR^{d\times d\times\VN\times \VN}$ is a block diagonal tensor with components 
\begin{align*}
\widetilde{\MBA}[\alp,\beta,\Vk,\Vl]=\del_{\Vk\Vl}A_{\alp\beta}(\x_\VN^\Vk);
\end{align*}
but one only needs to store the diagonals, which can be done in a  tensor  of shape $d\times d\times \VN$.

The numerical integration leads to an approximate formulation of the Galerkin approximation of \eqref{eq:WF}:
\begin{align*}
\text{find }u\in\cTN:\quad\bilfN{\nabla u_\VN}{\nabla v_\VN} = -\bilfN{\VE}{\nabla v_\VN},\quad\forall v_\VN\in\cTN;
\end{align*}
note that the approximation is exact for constant material coefficients $\T{A}$.
This formulation, that can be seen also as a collocation method \cite{ZeVoNoMa2010AFFTH}, is equivalent to the original Moulinec and Suquet formulation \cite{Moulinec1994FFT} in the sense that the solution vectors coincide \cite{VoZeMa2014FFTH}. 
However, the formulation here builds on the variational formulation \cite{VoZeMa2014FFTH} solved for the potential field (instead of gradient one).

The combination of numerical integration and differentiation of trigonometric polynomials \eqref{eq:grad_trig} allows to approximate the bilinear form in terms of the nodal values of potential fields
\begin{align*}
\bilfN{\nabla u_\VN}{\nabla v_\VN} 
=
\scal{\widetilde{\MBA}\iFFTd{\VN}\widehat{\nabla}_\VN\FFT{\VN}\MBu}{\iFFTd{\VN}\widehat{\nabla}_\VN\FFT{\VN}\MBv}_{\sR^{d\times\VN}}.
\end{align*}

In order to deduce the linear system, all operators acting on test vectors $\MBv_\VN$ are moved to the trial vector $\MBu_\VN$ as adjoint operators to reveal
the linear system in the original space 
\begin{align*}
\iFFT{\VN} \widehat{\nabla}_\VN^* \FFTd{\VN} \widetilde{\MBA} \iFFTd{\VN} \widehat{\nabla}_\VN \FFT{\VN} \MBu = -\iFFT{\VN} \widehat{\nabla}_\VN^* \FFTd{\VN} \widetilde{\MBA} \MBE,
\end{align*}
where $\MBE\in\sR^{d\times \VN}$ is constant with components $\MBE[\alp,\Vk]=E_\alp$.
One may notice that the system can be solved in Fourier space to save one computation of FFT and its inverse, which leads to the linear system in Fourier space
\begin{align}
\label{eq:ls_GaNi_fourier}
\widehat{\nabla}_\VN^* \FFTd{\VN} \widetilde{\MBA} \iFFTd{\VN} \widehat{\nabla}_\VN \hMBu = - \widehat{\nabla}_\VN^* \FFTd{\VN} \widetilde{\MBA} \MBE.
\end{align}


\subsubsection{Fourier-Galerkin method with exact integration (GA)}
\label{sec:FFTHga}
For many types of material coefficients \eqref{eq:A} and basis functions, there is a possibility to integrate the bilinear forms in the weak formulation exactly, which leads to a Galerkin approximation with exact integration
\begin{align}
\label{eq:Ga}
\text{find }u\in\cTN:\quad\bilf{\nabla u_\VN}{\nabla v_\VN} = \bilf{\VE}{\nabla v_\VN}\quad\forall v_\VN\in\cTN.
\end{align}
However, the exact integration of the Fourier-Galerkin formulation, in contrast to FEM, leads to a full linear system, which can be overcome with a double-grid integration with projection (DoGIP) \cite{VoZeMa2014GarBounds,Vondrejc2015FFTimproved}.
The DoGIP is a general method applicable also within the finite element method \cite{Vondrejc2017DoGIP-FEM}.
The original evaluation of the material law on a grid of size $\VN$ is reformulated as an evaluation on a double grid $2\VN-\V{1}$ with modified material coefficients; they can be expressed as a modification of the original material coefficients.

The main idea relies on expressing gradients of the trial and a test function together
\begin{align*}
\nabla u_\VN(\x)\otimes\nabla v_\VN(\x) = \Ve_\VN(\x)\otimes\Vw_\VN(\x) = \sum_{\Vk\in\ZtNd}{\MBe[:,\Vk]\otimes \MB{w}[:,\Vk] \varphi^{ \Vk }_\tVN(\x)}
\end{align*}
with respect to the basis of the double grid space consisting of trigonometric polynomials with doubled frequencies $\cT_\tVN$; the arrays $\MBe$ and $\MB{w}$ store the values of the trigonometric polynomials on the double grid, e.g.\ $\M{e}[\alp,\Vk] = \Ve_{\VN,\alp}(\x_\tVN^\Vk)$ for $\alp\in\{1,\dotsc,d\}$ and $\Vk\in\ZtNd$.
Then the bilinear form can be expressed on the double grid
\begin{align*}
\bilf{\Ve_\VN}{\Vw_\VN} 
&= \sum_{\Vk\in\ZtNd} \int_{\puc} \T{A}(\x) \varphi^{\Vk}_\tVN(\x)\D{\x}:\MBe_\tVN[:,\Vk]\otimes \MB{w}_\tVN[:,\Vk]
= \scal{\MBA\MBe}{\MB{w}}_{\sR^{d\times(\tVN)}}
\end{align*}
where $:$ is a double contraction between two matrices of size $d\times d$ and the material coefficients are defined as
\begin{align*}
\MA[\alp,\beta,\Vk,\Vl] = \del_{\Vk\Vl}\int_{\puc} A_{\alp\beta}(\x) \varphi^{\Vk}_\tVN(\x)\D{\x}\quad\text{for }\alp,\beta\in\{1,\dotsc,d\}\text{ and }\Vk,\Vl\in\ZtNd.
\end{align*}
This integration can be performed exactly for a large class of material coefficients.
In particular in \cite{Vondrejc2015FFTimproved,VoZeMa2014GarBounds}, square or circular inclusions have been considered, as well as image-based composites, materials with coefficients constant or bilinear over pixels (voxels in 3D).
Moreover, the evaluation of modified material coefficients can be performed effectively by FFT.

In order to derive the linear system, we have to still describe the interpolation from the original to the double grid space. As the spaces of trigonometric polynomials are nested $\cTN\subset\cT_{\V{M}}$ for $\VN<\V{M}$ (element-wise), we can just inject the polynomial to the bigger space by adding trigonometric polynomials with zero Fourier coefficients. This can be represented by the zero-padding injection operator
$\mathcal{I}:\sC^{d\times\VN}\rightarrow\sC^{d\times(\tVN)}$, defined as
\begin{align*}
(\mathcal{I}\hMBw)[:,\Vk] =
\begin{cases}
\hMBw[:,\Vk],&\text{for }\Vk\in\ZNd
\\
\V{0}&\text{for }\Vk\in\ZtNd\setminus\ZNd
\end{cases}.
\end{align*}
Its adjoint operator $\mathcal{I}^* :\sC^{d\times(\tVN)}\rightarrow\sC^{d\times\VN}$ just removes the frequencies $\Vk\in\ZtNd\setminus\ZNd$, i.e. projects on the $\Vk\in\ZNd$.

This allows us to deduce the linear system with exact integration
\begin{align}
\label{eq:ls_Ga_fourier}
\widehat{\nabla}_\VN^* \mathcal{I}^* \FFTd{\tVN} \MBA \iFFTd{\tVN} \mathcal{I} \widehat{\nabla}_\VN \hMBu = - \widehat{\nabla}_\VN^* \mathcal{I}^* \FFTd{\tVN} \MBA \MBE,
\end{align}
which has very similar structure compared to the scheme based on numerical integration \eqref{eq:ls_GaNi_fourier}.

\subsection{Preconditioning}
Following the recent paper \cite{LaPuVoZe2019precondFFT}, the preconditioning of both linear systems \eqref{eq:ls_GaNi_fourier} and \eqref{eq:ls_Ga_fourier} is based on a Laplacian expressed in the Fourier domain as
\begin{align*}
\hMBP[\Vk,\Vl] = \del_{\Vk\Vl}\Vk\cdot\Vl\quad\text{for }\Vk,\Vl\in\ZNd,
\end{align*}
which is a simple diagonal preconditioner. Its inverse is given by the Moore-Penrose pseudoinverse $\hMBP^{-1}[\Vk,\Vl]=\del_{\Vk\Vl}\frac{1}{\Vk\cdot\Vk}$ for $\Vk\in\ZNd\setminus\{\V{0}\}$ and $\hMBP^{-1}[\V{0},\V{0}]=0$; the latter condition enforces the zero-mean property of the approximated vectors.
The preconditioned systems are explicitly stated for both discretisation schemes
\begin{subequations}
\label{eq:Pls_fourier}
\begin{align}
\label{eq:Pls_GaNi_fourier}
\hMBP^{-1}\widehat{\nabla}_\VN^* \FFTd{\VN} \widetilde{\MBA} \iFFTd{\VN} \widehat{\nabla}_\VN \hMBu &= - \hMBP^{-1}\widehat{\nabla}_\VN^* \FFTd{\VN} \widetilde{\MBA} \MBE,
\end{align}
for the preconditioning of \eqref{eq:ls_GaNi_fourier}, and
\begin{align}
\label{eq:Pls_Ga_fourier}
\hMBP^{-1}\widehat{\nabla}_\VN^* \mathcal{I}^* \FFTd{\tVN} \MBA \iFFTd{\tVN} \mathcal{I} \widehat{\nabla}_\VN \hMBu &= - \hMBP^{-1} \widehat{\nabla}_\VN^* \mathcal{I}^* \FFTd{\tVN} \MBA \MBE.
\end{align}
for the preconditioning of \eqref{eq:ls_Ga_fourier}.
\end{subequations}

\section{FFT-based methods with low-rank approximations}\label{sec:fft-based-methods-with-low-rank-approximations}
Applying low-rank approximation techniques is of particular interest for problems with a huge number of degrees of freedom.
The low-rank approximations can not only furnish a posterior data compression of the solution array, but also reduce computational complexity by exploiting low-rank format representations in the solution process.
For the latter one needs some operations such as additions, element-wise multiplication, and the fast Fourier transform (FFT) to be implemented on tensors in low-rank format.
In this section we introduce an FFT-based solution process incorporating low-rank representations of tensors.
In the following sub-section~\ref{sec:overview-of-low-rank-formats}, the low-rank approximation formats are summarised along the corresponding operations;
details can be found in textbooks or in appendix~\ref{sec:low-rank-tensor-approximations}.
Then the application of low-rank approximation for the Fourier-Galerkin method is described and discussed in sub-section~\ref{sec:linear-systems-in-low-rank-formats}, and the suitable linear solvers in sub-section~\ref{sec:linear-solvers}.

\subsection{Overview of low-rank formats}\label{sec:overview-of-low-rank-formats}
Here we give a brief introduction of three types of low-rank tensors that are applied in this work, they are of canonical polyadic (CP), Tucker, and tensor train format respectively.
The CP format is only used  for tensors of order two because of  its intrinsic difficulty in finding optimal approximation  for tensors    with  higher   order.
The necessity and impact of rank truncation is also emphasized.
Interested readers are provided by more details about the operations on the low-rank tensors in the Appendix~\ref{sec:low-rank-tensor-approximations}.

\subsubsection{Canonical polyadic format}

A CP $r$-term approximation of a tensor $\MB v \in \mathbb K^{N_1 \times \cdots \times N_d}$ (the field $\mathbb K$ is $\mathbb R$ or $\mathbb C$) is a sum of $r$ rank-$1$ tensors. In this work the CP format is only used for tensors of order two ($d=2$), i.e. matrices.
In this case the representation has the form:
\begin{align*}
\MBv \approx \widetilde{\MBv}
&=\sum^r_{i=1} \M{c}[i] \bas{1}[i]\otimes\bas{2}[i],
\end{align*}
where $\MB{c}\in\sR^r$ stores the coefficients with respect to vectors $\bas{j} \in \sK^{r\times N_j}$ in the directions of indices $j$.
A low-rank representation for order-2 tensors (matrices) can be obtained by various matrix factorizing methods, among which the Singular Value Decomposition (SVD) is prominent as it provides a factorization that minimises the Frobenius-norm error of an $r$-term approximation.
The level of compression (reduction of memory requirements) depends on the rank $r$.
In order to find a solution in such a low-rank form, it requires to perform several operations occurring in the Fourier Galerkin method, particularly the FFT and element-wise multiplication.

The linearity and the tensor-product structure of the Fourier transform facilitates to express $d$-dimensional FFT of a tensor  (of order $d$)  as the sum of tensor products of $1$-dimensional FFTs, i.e.,
\begin{align*}
\FFT{\VN} (\widetilde{\MBv})
&= \sum^r_{i=1}  \M{c}[i] \FFT{N_1}(\bas{1}[i] )\otimes \FFT{N_2}(\bas{2}[i]).
\end{align*}
For the same number of tensor components in all directions $j$, i.e.\ 
$N_j=N$, 
this $d$-dimensional FFT algorithm has a complexity $O(drN\log N)$, which is much better than $O(dN^d\log N)$ for the full tensor, when the rank $r$ is kept low.
Note that this operation does not change the rank of a transformed tensor.

Another operation that occurs in the Fourier-Galerkin method is the sum and the element-wise (Hadamard) product of two tensors in low-rank format.
In the case of the CP format it is computed as:
\begin{align*}
\widetilde{\MBv} + \widetilde{\MBw} &= \sum^r_{i=1} \M{c}_{\MBv}[i] \left (\bas{1}_{\MBv}[i]\otimes \bas{2}_{\MBv}[i] \right )+  \sum^s_{k=1}  \M{c}_{\MBw}[k] \left (\bas{1}_{\MBw}[k]\otimes\bas{2}_{\MBw}[k]  \right ),
\\
\widetilde{\MBv} \odot \widetilde{\MBw} &= \sum^r_{i=1} \sum^s_{k=1} \M{c}_{\MBv}[i] \M{c}_{\MBw}[k] \left (\bas{1}_{\MBv}[i] \odot \bas{1}_{\MBw}[k] \right)\otimes \left (\bas{2}_{\MBv}[i] \odot \bas{2}_{\MBw}[k] \right).
\end{align*}
While addition of two tensor costs no floating point operations and only requires more memory, the element-wise multiplication has a complexity of $O(r s dN)$, which is significantly less than the $N^d$ operations for full tensors, especially when the ranks $r$ and $s$ are much smaller than $N$.

\subsubsection{Tucker format}
The decomposition of higher order tensors has many variants.
The Tucker format representation is linked to the definition of a tensor subspace $\sV=\bigotimes_{j=1}^d\sV^{j}$ where $\sV^j$ is a subspace of $\sR^{N_j}$ generated
by the span of vectors $\{\bas{j}[i]\,|\, i=1,\dotsc,r_j\}$; these vectors, which may be a frame, are typically chosen as an orthogonal or orthonormal basis. The Tucker format is then a linear combination of tensor products of all possible combinations of basis vectors in different directions, i.e.
\begin{align*}
\MBv \approx \sum^{r_1}_{i_1=1}\cdots\sum^{r_d}_{i_d=1} \M{c}[i_1,\dotsc,i_d] \bigotimes^d_{j=1}\bas{j}[i_j] \in \sR^\VN,
\end{align*}
where the core $\MB{c}\in\bigotimes_{\alp=1}^d\sR^{r_\alp}$ is a tensor of order $d$.
The CP format is then a special form of the Tucker format with a diagonal core.
Note that naturally there can be different number of basis vectors in different directions. 

\subsubsection{The Tensor train (TT) format}
The tensor train is another format which is suitable for the decomposition of higher order tensor. The idea is based on recursive decompositions done sequentially along the  \emph{tensor's} individual spatial  dimensions.
 For tensors of order $3$, the decomposition of the tensor of size $N\times N\times N$ is computed in two steps.
Using the standard SVD algorithm, the decomposition is first computed on the reshaped matrix of size $N\times N^2$.
It is followed by the decomposition of the reshaped right-singular vectors, i.e. of the matrix of size $N\times N$. The above recursive decomposition thus leads to
\begin{align}
\label{eq:TT}
\MBv &= 
\sum^{r_1}_{i_1=1}\sum^{r_2}_{i_2=1} \bas{1}[1,:,i_1] \otimes \bas{2}[i_1,:,i_2] \otimes \bas{3}[i_2,:,1],
\end{align}
where the vectors $\bas{j}[i_{j-1},:,i_j]\in\sR^{N_j}$ are vectors in direction $j$.
The tensor's components can be explicitly written for $\Vk=(k_1,k_2,k_3)$ as
\begin{align*}
\MBv[\Vk] = \sum^{r_1}_{i_1=1}\sum^{r_2}_{i_2=1} \basc{1}[1,k_1,i_1] \basc{2}[i_1,k_2,i_2] \basc{3}[i_2,k_3,1].
\end{align*}
For $d=2$ it is identical to the CP format.

The tensor train format in \eqref{eq:TT} is again expressed as a linear combination of rank one tensors, on which a $d$-dimensional FFT can be applied through a series of one-dimensional FFTs on the train \emph{carriages} along the second index, i.e.\ applied on the vectors $\bas{j}[i_{j-1},:,i_j]\in\mathbb{K}^{N_j}$ for all $i_j$.
The operations addition or element-wise multiplication are discussed in the Appendix~\ref{sec:tensor-train-format}.

\subsubsection{Rank truncation}
Rank truncation is the way to reduce computational complexity by a reasonable compromise in the precision of the low-rank approximations.
It is particularly necessitated by the fact that operations on low-rank tensors like addition and element-wise multiplication usually inflate the representation rank $\V r$, potentially at a very fast rate, which is detrimental to a fast computation. 
On the other hand, in the resulted representation, a large part of the $\V r$ terms are not essential and can be given up without or with minor loss of accuracy, if done correctly.

Rank truncations of tensors in the three low-rank formats are all based on QR decomposition, SVD, or high order SVD (HOSVD) \cite{Hackbusch2012book}, which provide optimal or suboptimal truncations and error estimates.

Other truncations are also possible. 
Particularly, the element-wise multiplication of two tensors with rank $r$ results in a tensor of rank $s=r^2$, which is truncated with computational complexity  $O(Ns^2)$ for CP and Tucker and $O(Ns^3)$ for TT format.
In case of higher rank $r$ of the original tensors, the truncation become computational bottleneck.
To speed up the basis orthogonalization procedure,
the  basis with relatively small norms can also be removed before the orthogonalization to trade accuracy for efficiency.
This is usually beneficial in an iterative solver.

We supplement a more detailed introduction to the truncation procedure in each low-rank format in the Appendix~\ref{sec:low-rank-tensor-approximations}.

\subsection{Applications of  low-rank approximations  on  the linear systems}
\label{sec:linear-systems-in-low-rank-formats}
Here, we discuss the application of low-rank formats on the linear systems \eqref{eq:Pls_fourier}, which are again stated here for the reader's convenience
\begin{subequations}
\label{eq:Pls}
\begin{align}
\overbrace{\hMBP^{-1}\widehat{\nabla}_\VN^* \FFTd{\VN} \widetilde{\MBA} \iFFTd{\VN} \widehat{\nabla}_\VN }^{\tMBC} \hMBu &= \overbrace{- \hMBP^{-1}\widehat{\nabla}_\VN^* \FFTd{\VN} \widetilde{\MBA} \MBE}^{\tMBb},
\\
\underbrace{\hMBP^{-1}\widehat{\nabla}_\VN^* \mathcal{I}^* \FFTd{\tVN} \MBA \iFFTd{\tVN} \mathcal{I} \widehat{\nabla}_\VN}_{\MBC} \hMBu &= \underbrace{- \hMBP^{-1} \widehat{\nabla}_\VN^* \mathcal{I}^* \FFTd{\tVN} \MBA \MBE}_{\MBb}.
\end{align}
\end{subequations}

The solution vector $\MBu$ (or their Fourier coefficients $\hMBu=\MBF\MBu$) stores the values of the trigonometric polynomial on the $d$-dimensional regular discretisation grid.
Therefore the solution vector can be naturally represented as a   tensor of order $d$, which allows a low-rank representation.
In order to avoid the computation of the full tensor and its decomposition, the low-rank tensor $\hMBu$ is computed by a suitable iterative solver introduced in \ref{sec:linear-solvers}. 
It requires to perform matrix vector multiplication for a low-rank tensor $\MBv$, which is approximated as
\begin{subequations}
 \label{eq:Pls_lowrank}
\begin{align}
\tMBC\MBv &\approx \Tr \hMBP^{-1} \Tr \widehat{\nabla}_\VN^* \FFTd{\VN} \Tr \widetilde{\MBA} \iFFTd{\VN} \widehat{\nabla}_\VN \MBv,
\\
\MBC\MBv &\approx \Tr \hMBP^{-1} \Tr \widehat{\nabla}_\VN^* \mathcal{I}^* \FFTd{\tVN} \Tr \MBA \iFFTd{\tVN} \mathcal{I} \widehat{\nabla}_\VN \MBv;
\end{align}
similarly, the right-hand side of the linear systems is approximated by a low-rank tensors
\begin{align}
\tMBb &= - \Tr \hMBP^{-1}\widehat{\nabla}_\VN^* \FFTd{\VN} \widetilde{\MBA} \MBE,
\\
\MBb &= - \Tr\hMBP^{-1} \widehat{\nabla}_\VN^* \mathcal{I}^* \FFTd{\tVN} \MBA \MBE.
\end{align}
\end{subequations}
These approximations involve several operations in low-rank formats such as differentiation, divergence, Fourier transform, and the truncation operator $\Tr$, which keeps the rank $\Vr$ at an affordable level.
The operations are tabulated in the Table~\ref{tab:operations} together with the corresponding implementations in low-rank format and their impact on the rank $\V{r}$.

\begin{table}[H]
 \centering
 \begin{tabular}{ |l| c|c |}
  \hline  
  Operation    & low-rank tensor implementation &  Rank $\V r$ \\[5pt]
  \hline  \hline 
  Differentiation (gradient) &  element-wise multiplication & remains unchanged  \\
  \hline
  Divergence &  element-wise multiplication and addition& is increased   \\ 
  \hline 
  Evaluation of material law &   element-wise multiplication & is increased    \\
  \hline
 $d$-dimensional FFT & series of 1D FFTs  & remains unchanged \\ 
  \hline 
   Preconditioning & element-wise multiplication  & is increased \\ 
  \hline 
 \end{tabular} 
 \caption{Operations and their implementations in low-rank formats}
 \label{tab:operations}
\end{table}

Since the material coefficients $\widetilde{\MBA},\MBA$ and the preconditioner $\MBP^{-1}$ are diagonal or block-diagonal for non-isotropic material coefficients, the related matrix-vector multiplications are implemented as element-wise multiplications, which inevitably inflates the representation rank of the tensors in low-rank format.
We apply a rank truncation after each multiplication to keep the computational complexity at a relatively low level, while maintaining reasonable accuracy in the solution.

The application of the gradient and divergence in Fourier space is also implemented as  element-wise multiplications. 
The differentiation operator for trigonometric polynomials is by nature a rank-1 tensor in the form
\begin{align*}
\widehat{\nabla}_\VN
&=
[2\pi\imu \V{K}_1\otimes \V{1} \otimes \V{1},\V{1} \otimes 2\pi\imu \V{K}_2\otimes \V{1},\V{1} \otimes \V{1}\otimes 2\pi\imu \V{K}_3]
\end{align*}
in the 3D setting, where $\V{K}_\alp=(k\in \sZ\,;\, |k|<N/2)$ is a vector of all discrete frequencies in direction $\alp$.
So the corresponding element-wise multiplication keeps the rank of tensors unchanged.
However, for the divergence the contraction along the first component of $\widehat{\nabla}_\VN$ is provided by the operation addition of two low-rank formats, which increases the rank, and hence a truncation has to be performed.

The last operation that occurs in the system is the $d$-dimensional fast Fourier transform (FFT) which is efficiently evaluated using $1$-dimensional FFTs.
Moreover the rank of the tensor remains the same in this operation.

\subsection{Linear solvers}
\label{sec:linear-solvers}

For the full solver we have used preconditioned conjugate gradients, which is considered to be the best available solver for FFT-based homogenisation \cite{MiVoZe2016jcp,LaPuVoZe2019precondFFT}.
However, the linear systems with low-rank approximations require solvers that are insensitive to small perturbations, as the matrix-vector product is computed only approximately, due to the truncation of tensors.
Therefore conjugate gradient method that builds on the orthogonalisation of Krylov subspace vectors using a short-term recurrence relation is inappropriate.

The systems with low-rank approximations are solved here with minimal residual iteration \cite{Saad2003IMSL} which is closely related to Richardson iteration.
The latter is well established in the FFT-based community, as it corresponds to the original Moulinec-Suquet algorithm.
Both methods solve the linear system $\MB{C}\MBu = \MB{d}$, see \eqref{eq:Pls} and \eqref{eq:Pls_lowrank} for details, by the iteration
\begin{align*}
\MBu_{(i+1)} = \MBu_{(i)} + \omega \underbrace{(\MB{d} - \MB{C}\MBu_{(i)})}_\textrm{$\MBr_{(i)}$} = \bigl(\MB{I} - \omega \MB{C}\bigr)\MBu_{(i)} + \omega\MB{d}.
\end{align*}
In the Richardson iteration, the parameter $\omega$ is chosen such that the iteration matrix $(\MB{I} - \omega \MB{C})$ has a norm smaller than one to guarantee convergence. A fixed value $\omega$ is set on the basis of a priori knowledge about the extreme eigenvalues of the system matrix $\MB{C}$, i.e.\
\begin{align*}
\omega=\frac{2}{\lambda_{\min}(\MB{C})+\lambda_{\max}(\MB{C})},
\end{align*}
because it satisfies the minimal norm of the iterative matrix as proposed in \cite{Moulinec1998NMC} for FFT-based homogenisation.
Here $\lambda_{\min}(\MB{C})$ denotes the smallest positive eigenvalue, as the system matrix is only positive semidefinite.
In particular, the linear systems in \eqref{eq:Pls_fourier} contains one zero eigenvalue corresponding to the constant fields, while the linear systems that are formulated in traditional FFT-based homogenisation for gradients fields contain many zero eigenvalues corresponding to the eigenspace composed of divergence-free fields.
In both cases the solver produces the solution in the space of compatible fields.

In the minimal residual iteration, the parameter $\omega$ is chosen at each iteration as the minimizer of the next residual $\MBr_{(i+1)}$ over all increments of $\MBu$ in the direction of $\MBr_{(i)}$., i.e.\
\begin{align*}
\MBu_{(i+1)} = \MBu_{(i)} + \omega_{(i)} \MBr_{(i)},\;\;\;\; \mbox{with}\;\;
\omega_{(i)}=\frac{\left(\MB{C}\,\MBr_{(i)}, \MBr_{(i)} \right )}{\|\MB{C}\,\MBr_{(i)} \|^2}.
\end{align*}
We adopt the latter method in this work, because of our observation that the minimal residual iteration is more robust than Richardson iteration, for which we have observed a divergence when a massive truncation has been used during the iterations.
For a low-rank approximation of a solution vector, note that the solver has to deal with the matrix vector product $\MB{C}\MBu_{(i)}$, which is computed only approximately \eqref{eq:Pls_lowrank} to limit the growth of the solution rank.
The rank also grows by the operation addition during the iteration. Therefore, a truncation is included at each step of the low-rank variant of the minimal residual iteration, i.e.
\begin{align*}
\MBu_{(i+1)} = \Tr[\MBu_{(i)} + \omega_{(i)}(\MB{d} - \MB{C}\MBu_{(i)}) ].
\end{align*}

%
%

\section{Numerical results}\label{sec:numerical-results}
 The methodology described in the previous sections is  tested on several numerical examples with material parameters defined in section~\ref{sec:material-parameters}.
We compare two numerical homogenisation schemes: the Fourier-Galerkin method with numerical integration (GaNi) and a version with exact integration (Ga), described in sections \ref{sec:FFTHga} and \ref{sec:FFTHgani}.
The preconditioned linear systems stated in \eqref{eq:Pls_fourier} are solved by conjugate gradient method.
The same systems that are equipped with low-rank tensor approximations are solved by the minimal residual iteration, which is discussed in section~\ref{sec:linear-solvers}.

The numerical results were calculated using software FFTHomPy (FFT-based Homogenisation in Python), which is freely available at \href{https://github.com/vondrejc/FFTHomPy}{https://github.com/vondrejc/FFTHomPy}; the software contains examples, which are described in the following sections.

\subsection{Material parameters}
\label{sec:material-parameters}

\begin{figure}[htb]
 \centering
 \subfloat[Square inclusion (\squar) in 2D]{\includegraphics[scale=0.48]{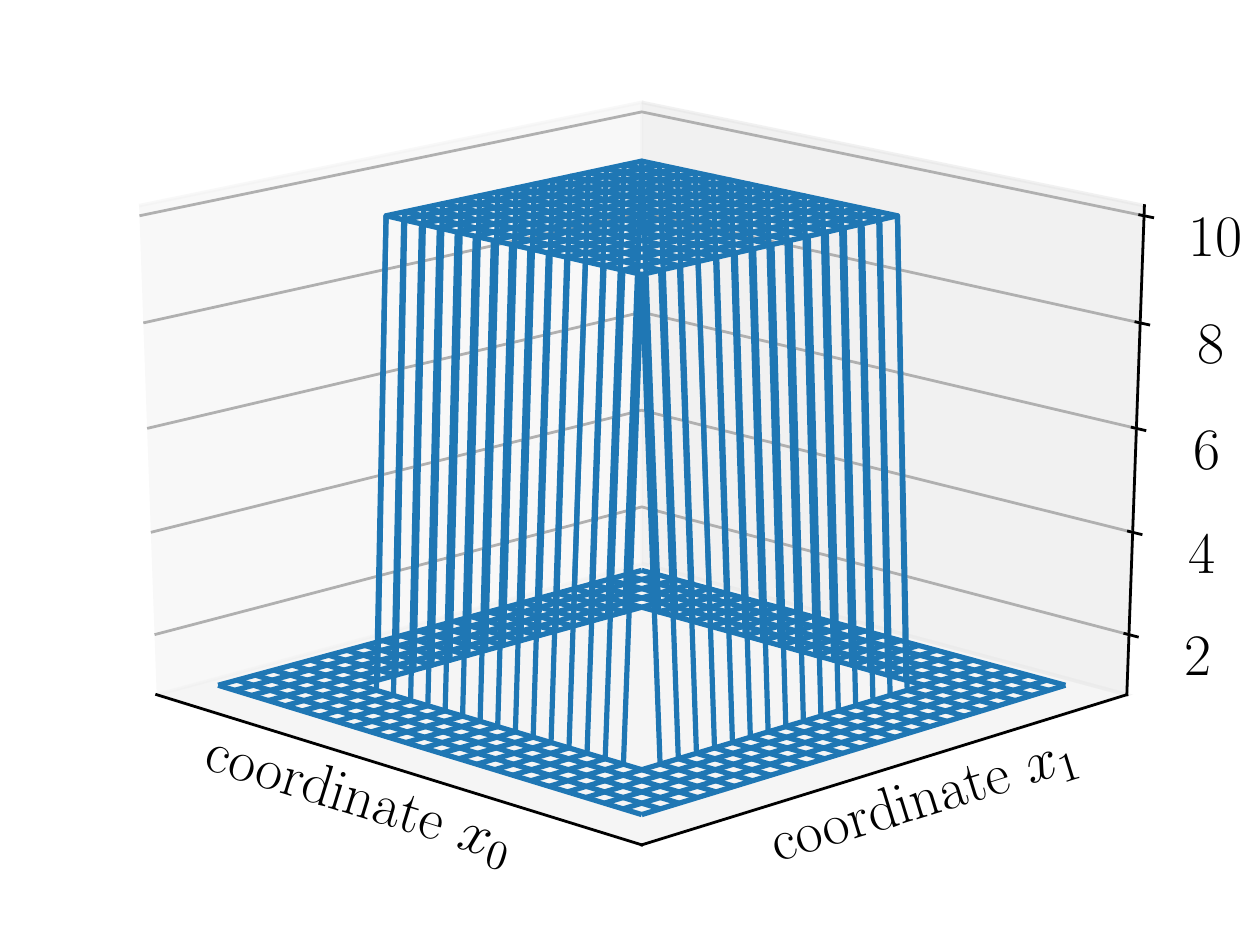}}
 \subfloat[Square (cube) inclusion (\squar) in 3D]{\includegraphics[scale=0.48]{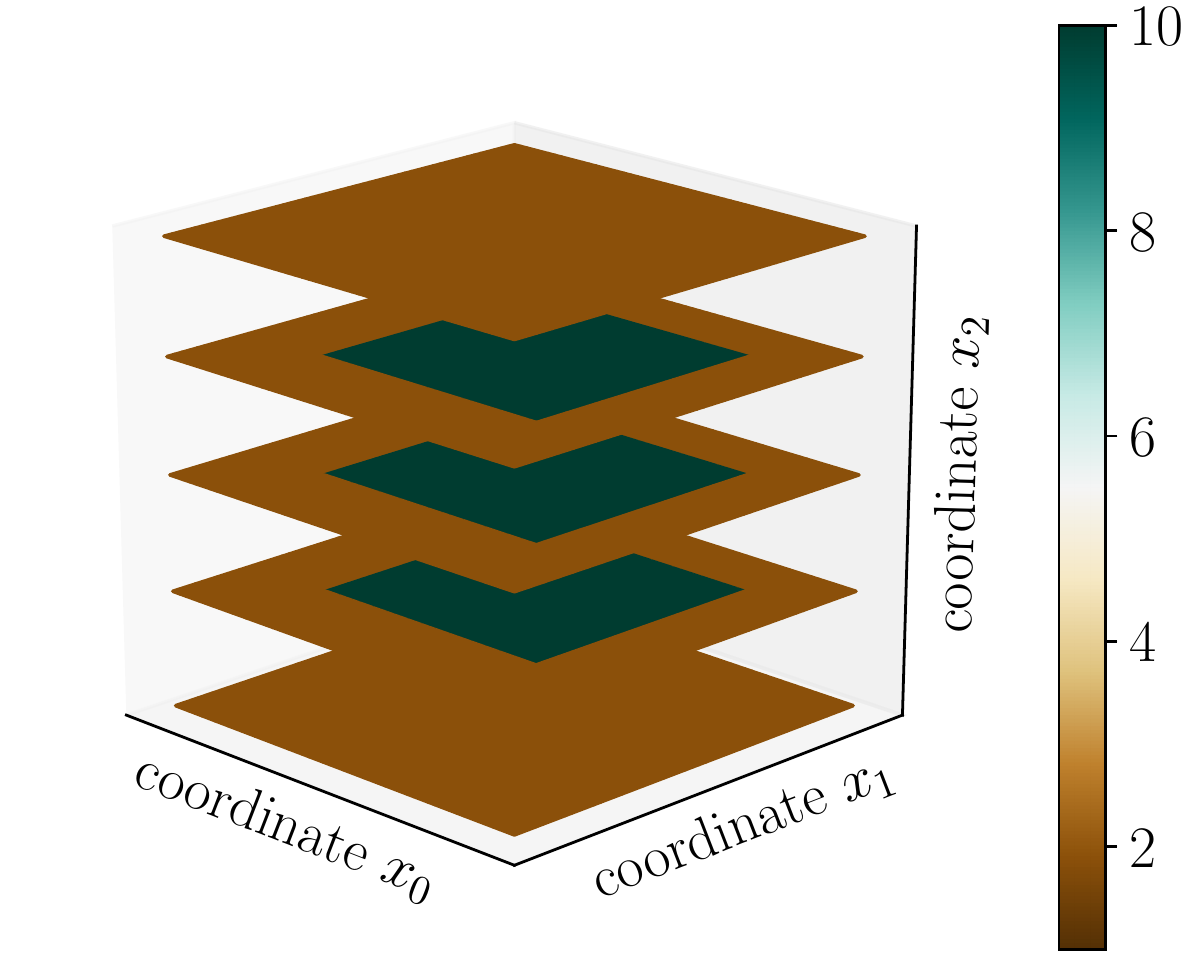}}
 \caption{Material coefficients \eqref{eq:Adef} of the square and the cube inclusion defined by \eqref{eq:shape_fun}.}
 \label{fig:shape_fun}
\end{figure}

\begin{figure}[htb]
 \centering
 \subfloat[Stochastic material (\stoch) in 2D]{\includegraphics[scale=0.48]{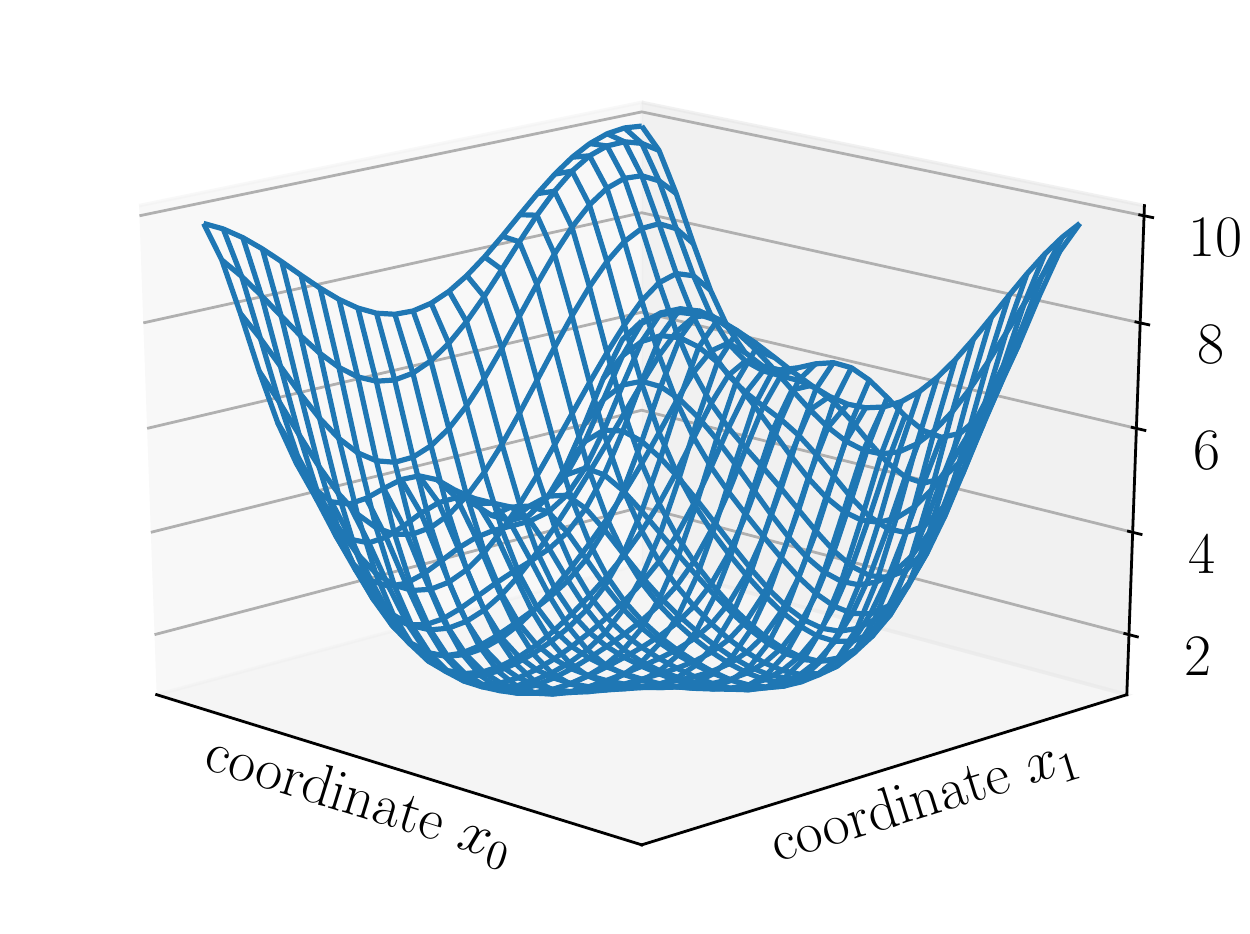}}
\subfloat[Stochastic material (\stoch) in 3D]{\includegraphics[scale=0.48]{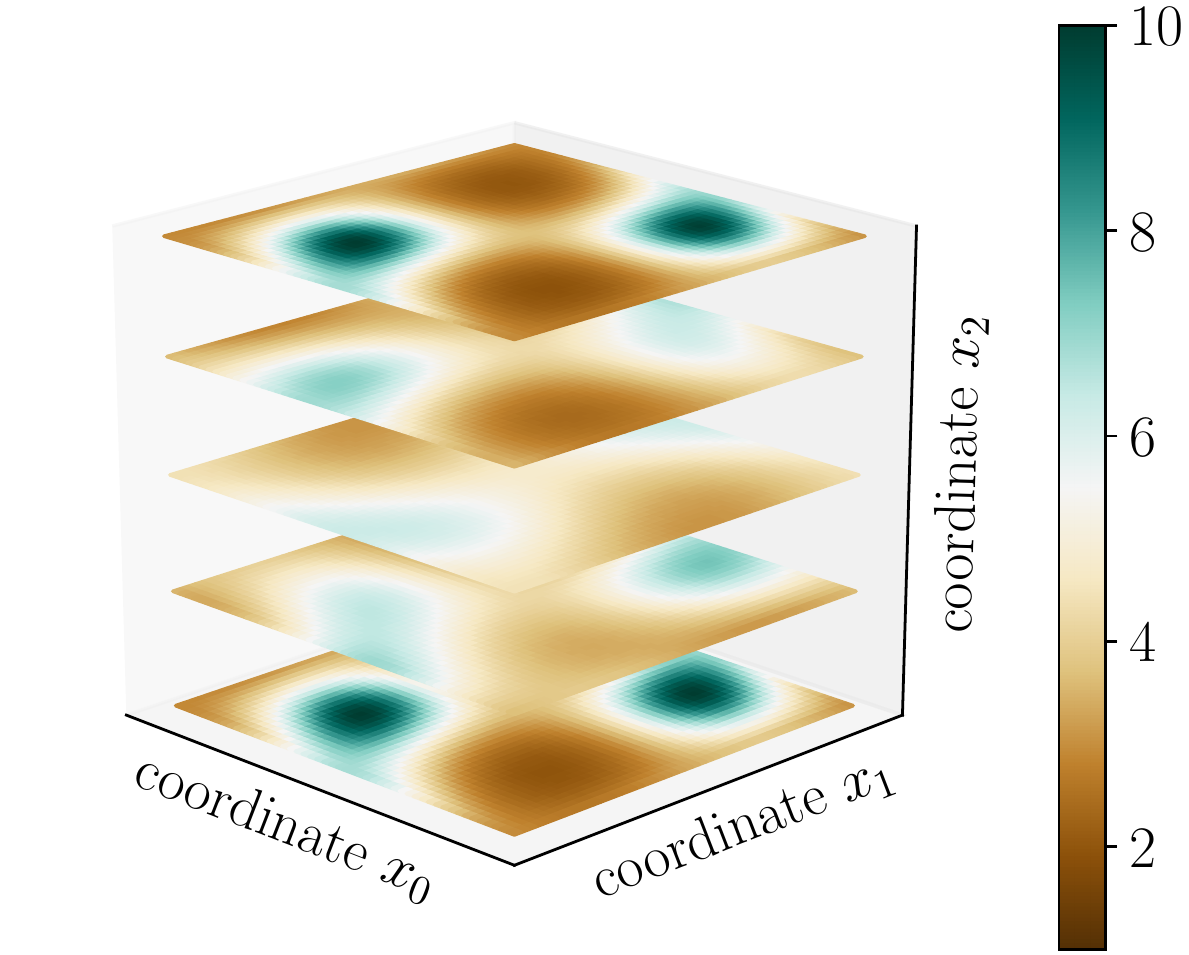}} 
 \caption{One sample of the stochastic material defined by \eqref{eq:shape_fun2}.}
 \label{fig:shape_fun2}
\end{figure}

Here, we present two material examples on which we did numerical tests. The first is defined as
\begin{align}
\label{eq:Adef}
\T{A}_{\squar}(\x) &= \T{I}(1+\rho \chi(\x))
\end{align}
where $\V{I}\in\sR^{d\times d}$ is the identity matrix and the parameter $\rho=10$ corresponds to a material contrast. The function $\chi:\puc\rightarrow\sR^d$ describing the topology of the inclusions is defined on a unit cell $\puc=(-\frac{1}{2},\frac{1}{2})^d$ as
\begin{align}
\label{eq:shape_fun}
\chi(\x) &=
\begin{cases}
1&\text{for }\x \text{ such that }x_i<0.3\text{ for }i=1,\dotsc,d,
\\
0&\text{otherwise}
\end{cases},
\end{align}
which is also depicted in 2D in Figure~\ref{fig:shape_fun}. The corresponding low-rank approximations have rank $2$ for all three formats (CP, Tucker, tensor train).

 As a second example, one sample of a stochastic material has been obtained using the  truncated Karhunen-Lo\`{e}ve expansion \cite{Adler2009} of the squared exponential Matérn covariance function \cite{Matern1986}.
 In order to obtain positive definite material coefficients, the exponential function has been applied on the expansion, which leads to the following form
\begin{align}
\label{eq:shape_fun2}
\T{A}_{\stoch}(\x) &=	\T{I}\exp\bigl(C + D\sum_{\Vk\in I} c[k] \varphi^\Vk(\x)\bigr).
\end{align}
The most important modes of the expansion has been selected ($20$ modes in $2D$ and $26$ modes in 3D) and the corresponding frequencies are collected in the index set $I$.
The coefficients $c[k]$ for $k\in I$ has been sampled from uniform distribution on the interval $[-0.5, 0.5]$. The constants $C$ and $D$ scales the material coefficients such that the minimal eigenvalue of $\T{A}$ is $1$, and the maximal $10$.
The particular sample that is used for the computation is plotted in Figure~\ref{fig:shape_fun2}.
The material coefficients were approximated in low-rank formats
with a rank set to $10$. 
For a comparison to the full solution, the full material coefficients have been recovered in order to compute exactly the same problem.

All the numerical problems have been computed with the same number of discretisation grids in each direction $\VN=[N,\dotsc,N]\in\sR^d$.

\subsection{Behaviour of linear systems during iterations}

\begin{figure}[htb]
 \centering
 \subfloat[\ga, CP, 2D, \squar]{\includegraphics[width=0.45\linewidth]{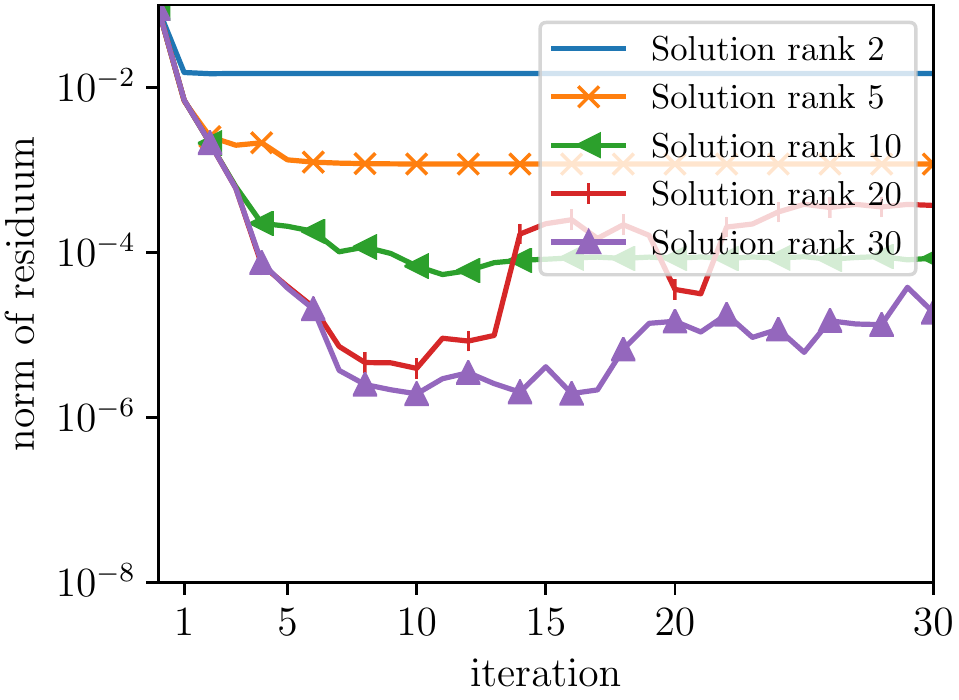}}
 \subfloat[\gani, CP, 2D, \stoch]{\includegraphics[width=0.45\linewidth]{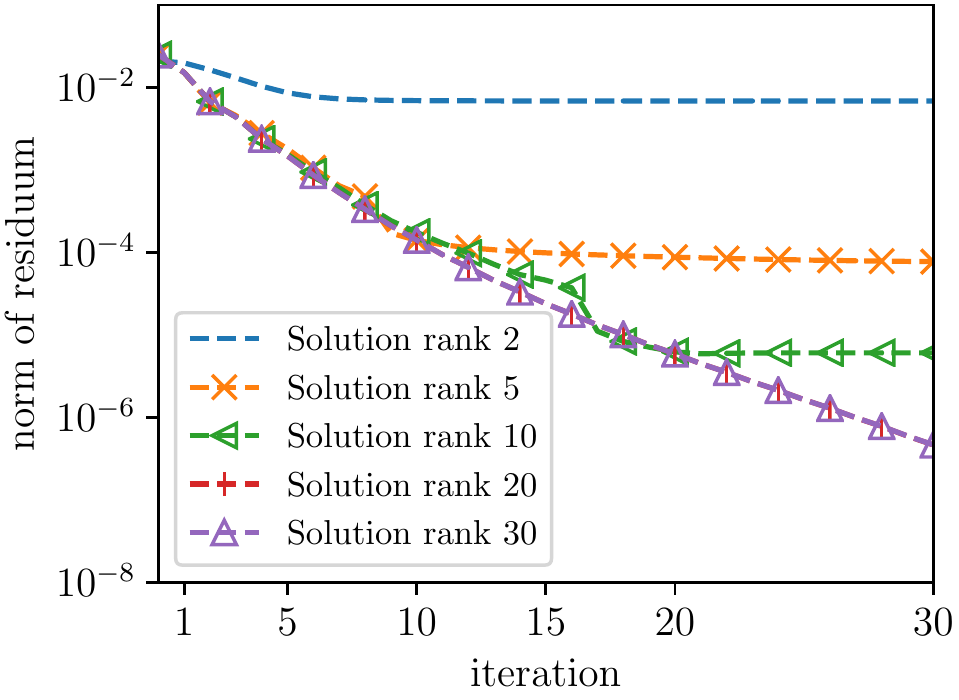}}
 \\
 \subfloat[\ga, \tucker, 3D, \squar]{\includegraphics[width=0.45\linewidth]{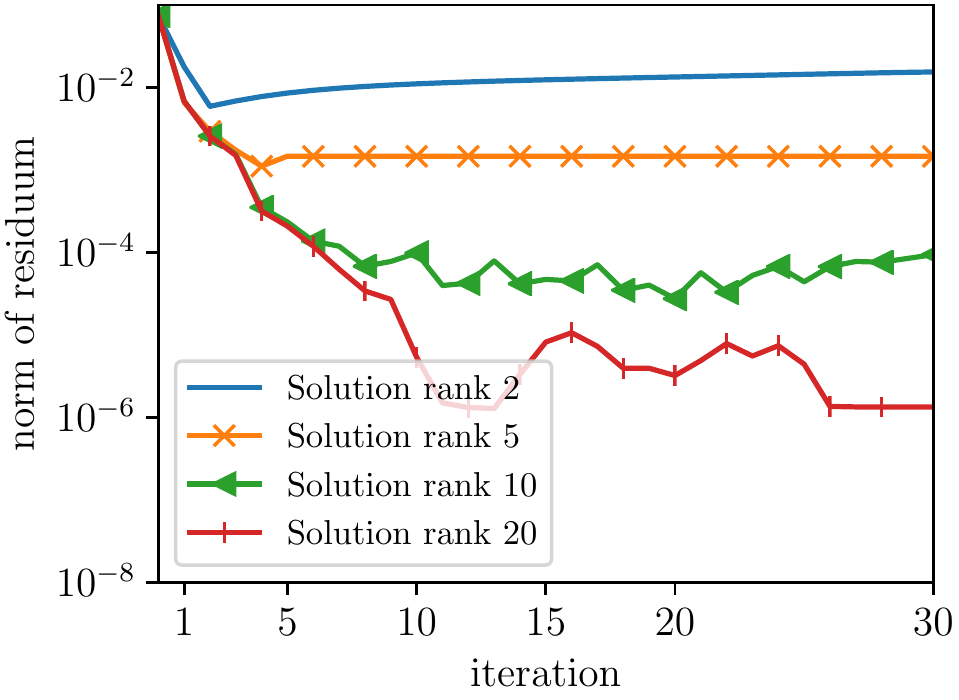}}
 \subfloat[\gani, \TT, 3D, \stoch]{\includegraphics[width=0.45\linewidth]{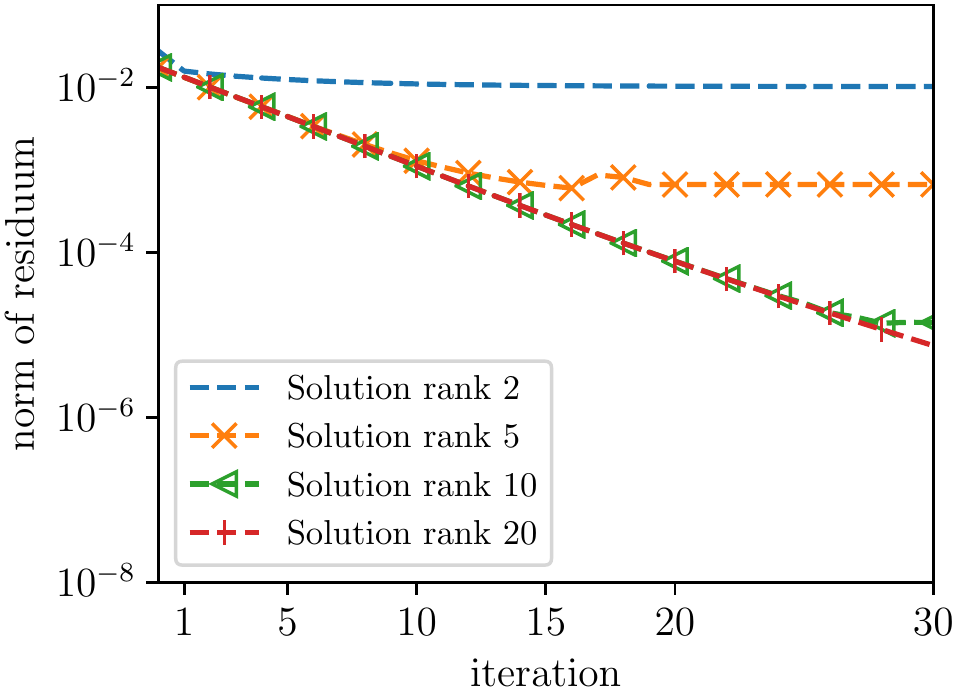}}
 \caption{Evolution of the norm of residua during minimal residuum iteration; computed in 2D for $N=1215$ and in 3D for $N=135$.}
 \label{fig:Residuum}
\end{figure}
The evolution of the norm during the minimal residual iteration is investigated because it describes well the character of the low-rank approximations.
The numerical results in Figure~\ref{fig:Residuum} depict the Euclidean norm of the residuum $\MB{r} = \MB{d}-\MB{C}\MBu\iter{i}$
\begin{align*}
\|\MB{r}\| = \Bigl(\sum_{\Vk\in\ZNd} \bigl|\M{r}[\Vk]\bigr|^2\Bigr)^{\frac{1}{2}}
\end{align*}
because it corresponds to the $L^2$-norm of the corresponding trigonometric polynomial.
Note that since the problem is solved in Fourier space, the residuum components agree with the Fourier coefficients of the corresponding trigonometric polynomial.

Although the truncation of the growing tensor's rank can be provided by a tolerance to an approximation error, it is difficult to set up the parameters properly during the solver.
Particularly it may happen that the rank significantly increase resulting in unnecessary computational demands, especially when the tensors are far away from the solution.
Therefore the truncation has been performed to a fixed rank.
The solution which is from a large dimensional space $\sR^\VN$  with the dimension $\prod_{\alp=1}^d N_\alp$  is approximated with a significantly smaller number of parameters.
Therefore there is always a residual error which can be diminished only by an increasing rank of the low-rank formats.
Note that the rank-one tensors occurring in all three low-rank formats are automatically computed by a solver and are thus suboptimal global basis vectors for the particular problem.
Therefore the method can be seen as a model order reduction technique. 

From the results in Figure~\ref{fig:Residuum}, 
we can observe that solutions with higher rank have larger potential in reducing the norm of residuum regardless the discretisation method ({\ga} and \gani), material problem ({\squar} and \stoch), or the low-rank format (CP, Tucker, TT).
This proposes a rank adapting solver that starts with a lower solution rank and increases the rank during the iterations.
We also notice that the norms of residuum during iterations decrease with higher rate for the problem with the square inclusion (material $\square$), however, the rate is more stable for the material S.
Although, the material $\square$ was systematically computed with GaNi method and material S with Ga, which is in accordance with the recommendation in \cite{VoGe2017FFTHvsFEM}, the discretisation method has no influence on the character of the behaviour during iterations.
These finding are in agreement with \cite{Matthies2012} analysing the stochastic linear systems and solvers approximated with low-rank approximations.

Note that the computation of the Frobenius norm of tensors in Tucker format is computationally demanding. Therefore, we have used the equivalent Frobenius norm of the Tucker's core, which can be computed much faster.

\subsection{Algebraic error of the low-rank approximations}
\begin{figure}[htb]
 \centering
 \subfloat[2D, \squar]{\includegraphics[width=0.48\linewidth]{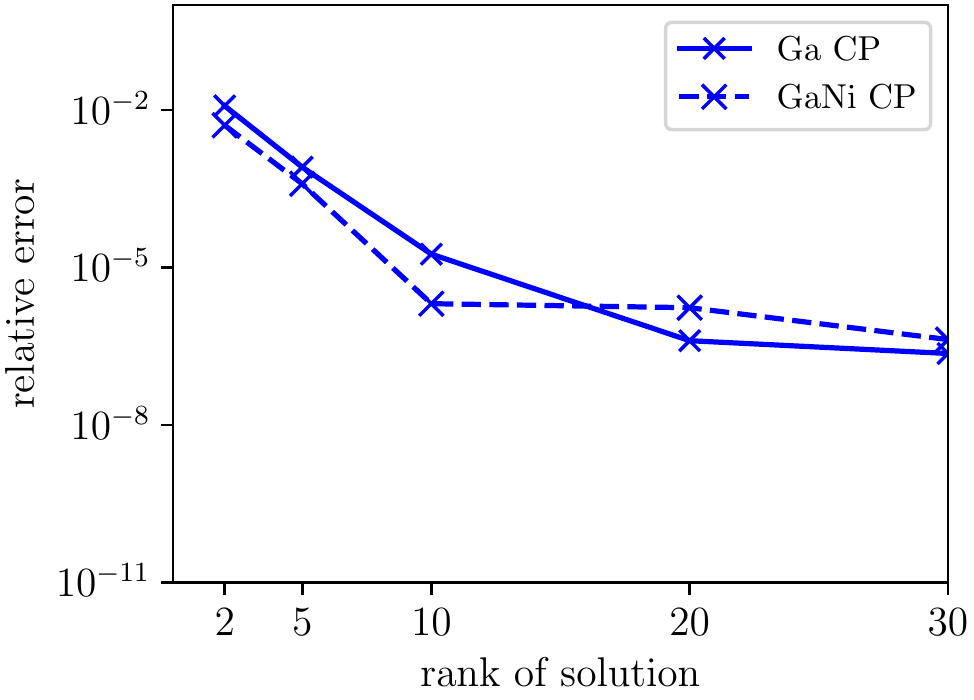}}
 \subfloat[2D, \stoch]{\includegraphics[width=0.49\linewidth]{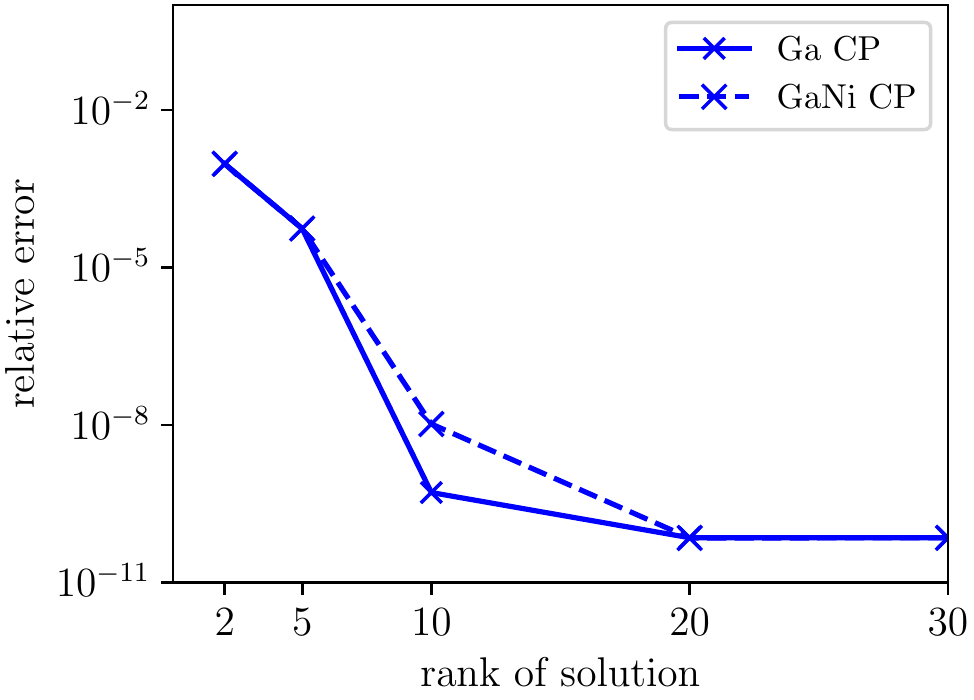}}
 \\
 \subfloat[3D, \squar]{\includegraphics[width=0.48\linewidth]{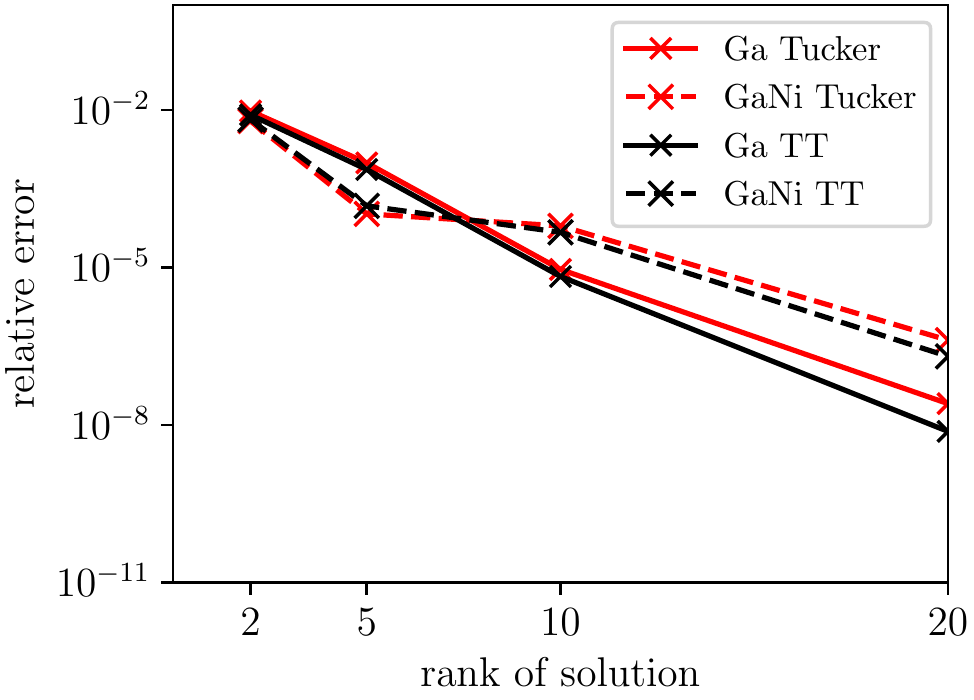}}
 \subfloat[3D, \stoch]{\includegraphics[width=0.49\linewidth]{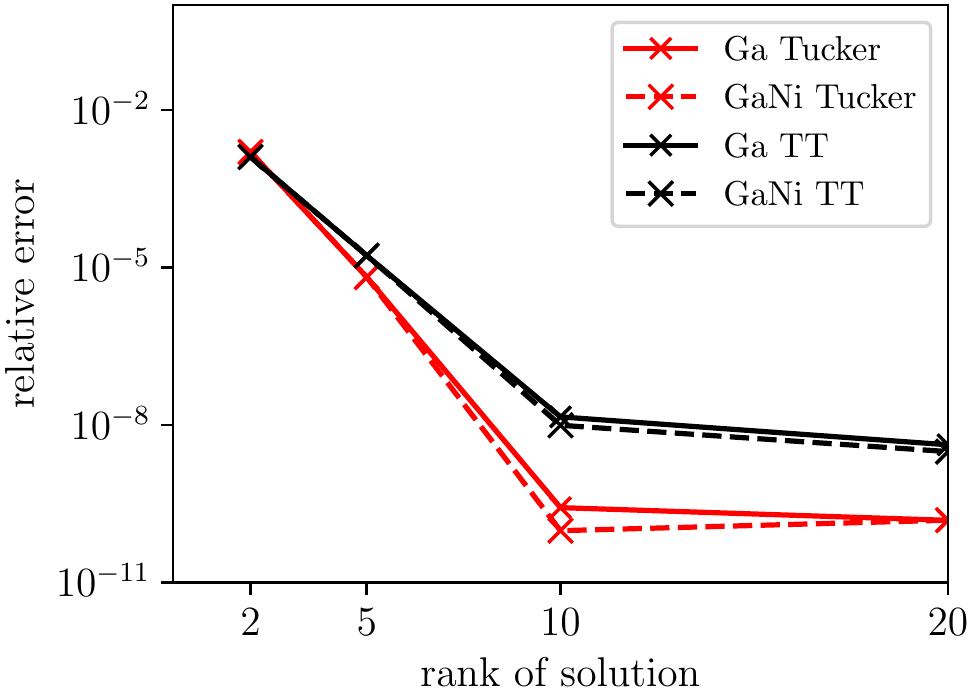}}
 \caption{Relative errors \eqref{eq:rel_error} of low-rank solutions computed in 2D for $N=1215$ and in 3D for $N=135$.}
 \label{fig:error}
\end{figure}

In the Figure~\ref{fig:error}, the approximation properties of the low-rank formats are depicted.
As an criterion, the relative algebraic error between the homogenised properties of low-rank solution $A_{\eff,\VN,r}$ and of the full solution $A_{\eff,\VN}$ has been used, i.e.
\begin{align}
\label{eq:rel_error}
\text{relative error}= \frac{A_{\eff,\VN}-A_{\eff,\VN,r}}{A_{\eff,\VN}}.
\end{align}
This is chosen because the error in the homogenised properties corresponds to the square of the energetic semi-norm (norm on zero-mean fields) of the algebraic error between the full solution and the low-rank approximation
\begin{align*}
\|u_\VN-u_{\VN,r}\|_A^2 = \bilf{\nabla u_\VN-\nabla u_{\VN,r}}{\nabla u_\VN-\nabla u_{\VN,r}} = A_{\eff,\VN,r}-A_{\eff,\VN};
\end{align*}
for the derivation see \cite[Appendix~D]{VoGe2017FFTHvsFEM}.
We also note that the full solution $u_\VN$ has been computed using conjugate gradients with high accuracy (tolerance $10^{-8}$ on the norm of the residuum) to obtain a solution that is close to the exact one.
 The low-rank solution has been obtained from minimal residual iteration, which was stopped when the residuum failed to be decreased.  The minimal residual iteration was used to provide low-rank solution with the minimal norm of residuum.  

We can observe that the results are again similar regardless of the discretisation method ({\ga} and \gani), material problem ({\squar} and \stoch), or the low-rank format (CP, Tucker, TT).
An increase in the solution rank leads to a significant reduction of the relative error.  However, the low-rank approximations of the material S reach the threshold error corresponding to the full approximate solution obtained from the conjugate gradients. It also shows that the low-rank method is more accurate for a problem with continuous material property (material \stoch) than for the one with discontinuous coefficients (material \squar).

\subsection{Memory and computational efficiencies}
 
\begin{table}[htb]
\centering
 \begin{tabular}{|l|lll|}
  \hline
  \backslashbox{Formats}{Operations} & Element-wise product & FFT$_d$ & Truncation\\ 
  \hline
  full &  $N^d$ & $\C O(N^d\log N$ )& --- \\
  CP  & $dNrs$  & $ \C O(dNr\log N)$ &$\C O(dNr^2)$ \\ 
  Tucker &   $dNrs+r^ds^d$ & $ \C O(dNr\log N) $ &$\C O(dNr^2 + r^{d+1})$ \\ 
  TT  & $dNr^2s^2$  & $ \C O(dNr^2\log N)$& $\C O(dNr^3)$   \\
  \hline
 \end{tabular}
 \caption{Asymptotic computational complexities in terms of floating point multiplications. The operations are performed on full tensors of order $d$ and shape $(N,\dotsc,N)$, and on the same tensors in their CP, Tucker, and tensor-train (TT) formats with maximum rank $r$ and $s$ ($s$ for the second operand in a binary operation).}
 \label{tab:complexities}
\end{table}

\begin{table}[htb]
\centering
 \begin{tabular}{|l|c|}
  \hline
  format & memory requirements\\ 
  \hline
  full & $N^d$ \\
  CP & $dNr$ \\ 
  Tucker & $dNr+r^d$ \\ 
  TT &$2Nr+(d-2)Nr^2$ \\
  \hline
 \end{tabular}
\caption{Memory requirements to store tensors of order $d$ with shape $(N,\dotsc,N)$ for full, CP, Tucker, and tensor-train (TT) formats with maximum rank $r$.}
\label{tab:memory}
\end{table}

Here, we discuss the computational and memory requirements to resolve the linear system using low-rank approximations.
Additionally  to the previous examples, the CPU times and approximation properties of low-rank formats were tested for an anisotropic material. The heterogeneous material coefficients $\TA_\squar$ and $\TA_\stoch$ were modified by adding a spatially constant anisotropic material tensor $\T{B}$, i.e.
 \begin{align*}
 \widetilde{\TA}_{\bullet}(\x) = \TA_\bullet(\x) + \T{B},
 \end{align*}
 where the matrices
 \begin{align*}
 \T{B} &=
 \begin{pmatrix}
 5.5 & -4.5 \\
 -4.5 & 5.5 
 \end{pmatrix},
 &
 \T{B} &=
 \begin{pmatrix}
 4.25    &   -3.25   &    -1.25\sqrt{2}\\
 -3.25   &     4.25    &    1.25\sqrt{2}\\
 -1.25\sqrt{2} & 1.25\sqrt{2} & 7.5
 \end{pmatrix}
 \end{align*}
have eigenvalues $(1,10)$ in 2D and $(1,5,10)$ in 3D.

 \begin{figure}[htb]
\centering
\subfloat[2D, \squar, isotropic]{\includegraphics[width=0.49\linewidth]{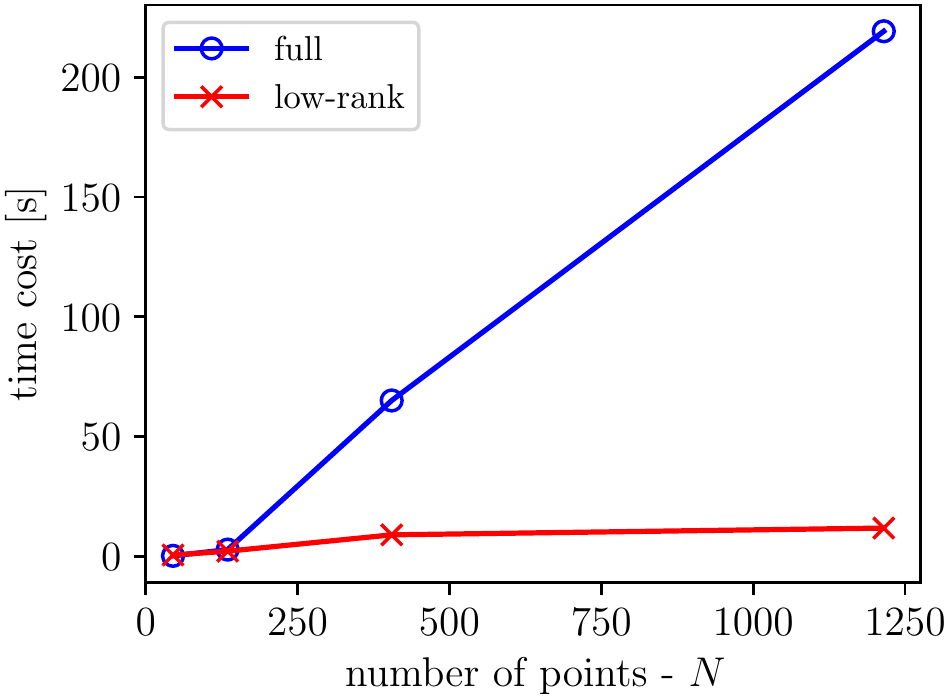}}
\subfloat[3D, \squar, isotropic]{\includegraphics[width=0.49\linewidth]{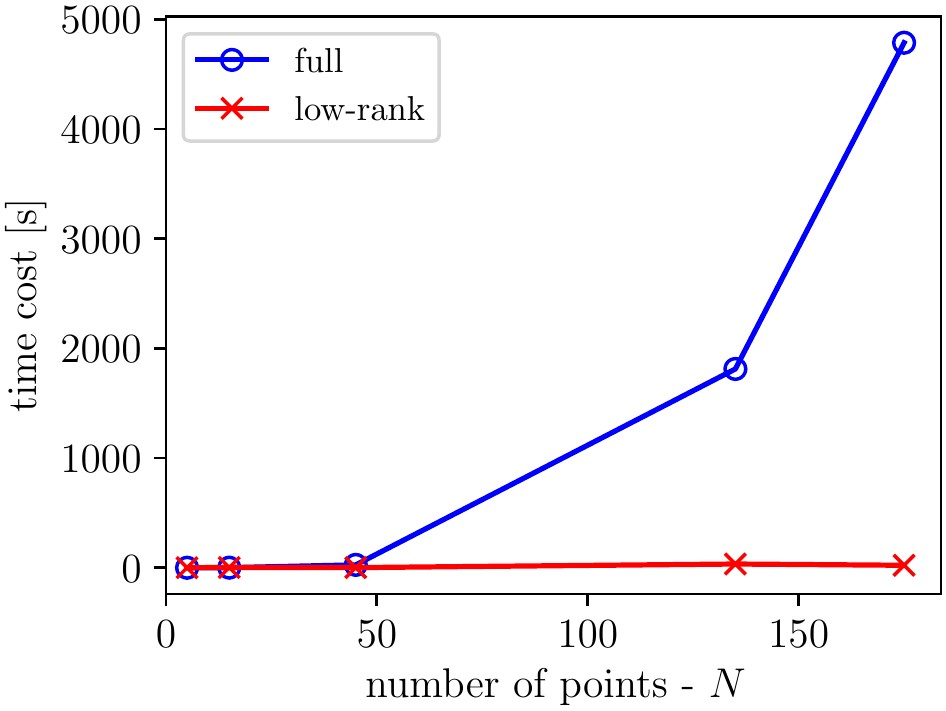}} \newline
\subfloat[2D, \squar, anisotropic]{\includegraphics[width=0.49\linewidth]{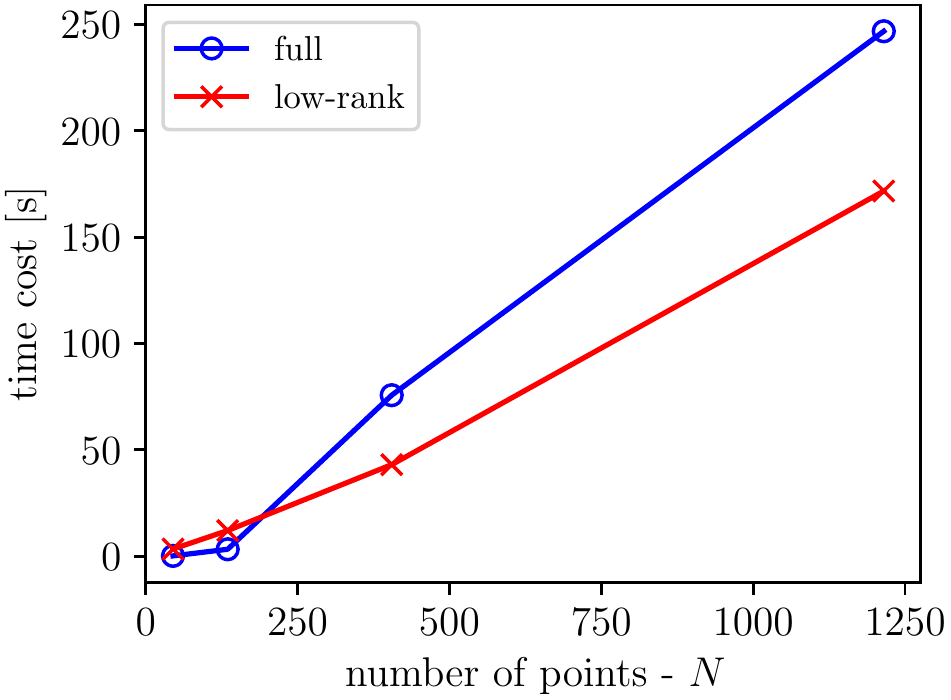}}
\subfloat[3D, \squar, anisotropic]{\includegraphics[width=0.49\linewidth]{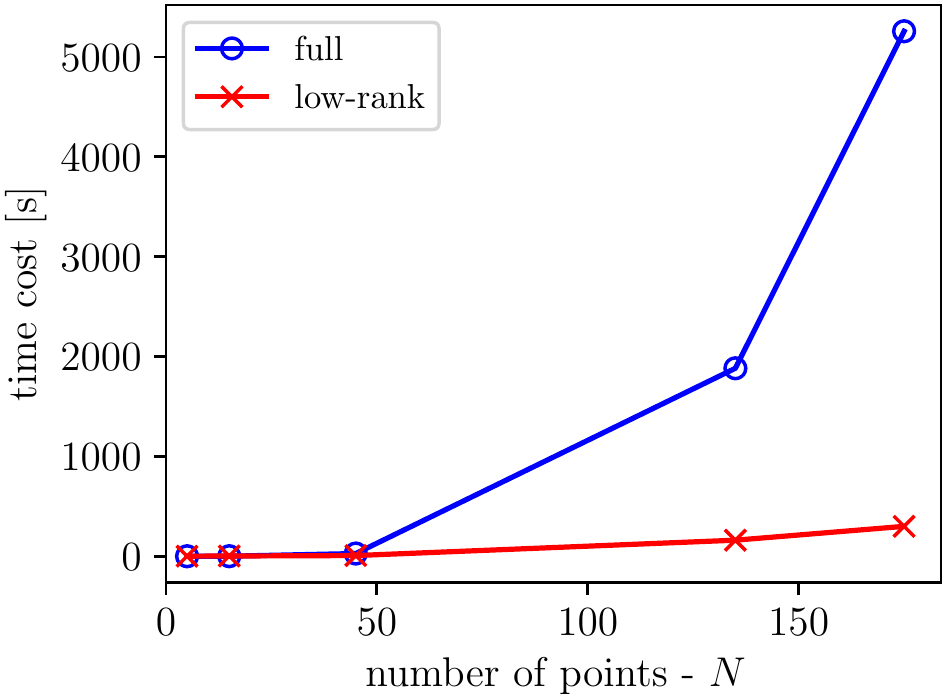}}
\caption{ The CPU time of Ga solver  to solve the problem with  isotropic (1st row) and anisotropic (2nd row) material \squar. The full solution has been computed on a grid of size $(N,\dotsc,N)$ while the low-rank solution on the grid $(3N,\dotsc,3N)$ with  various  solution ranks to achieve the same level of accuracy as the full scheme.  The stopping criterion of conjugate gradients for the full solver was set to $10^{-6}$ on the norm of residuum.} 
\label{fig:time_material_0_3}
\end{figure}
 \begin{figure}[htb]
\centering
\subfloat[2D, S, isotropic]{\includegraphics[width=0.49\linewidth]{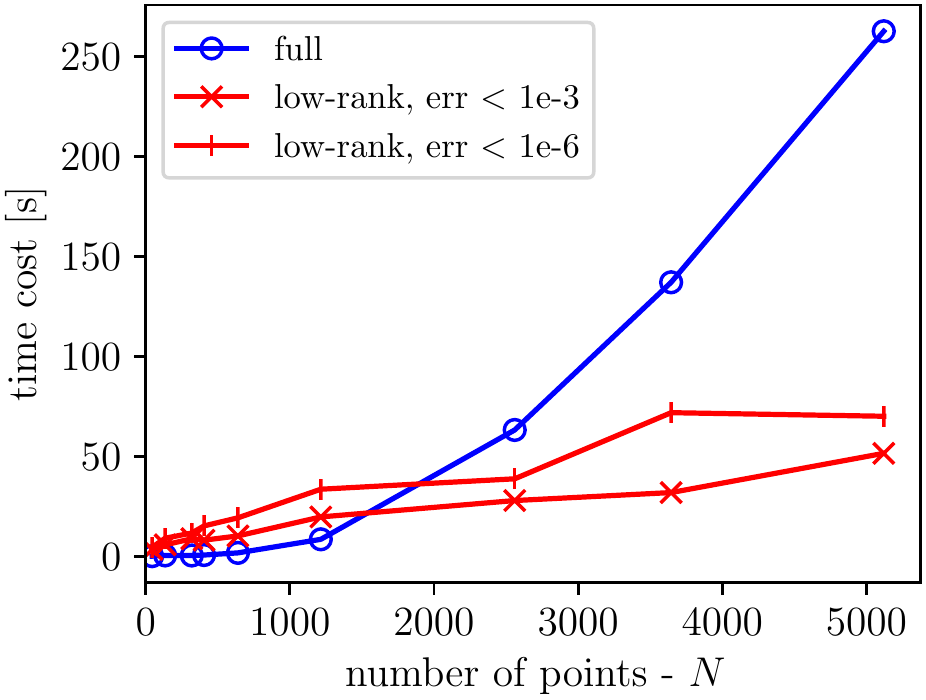}}
\subfloat[3D, S, isotropic]{\includegraphics[width=0.49\linewidth]{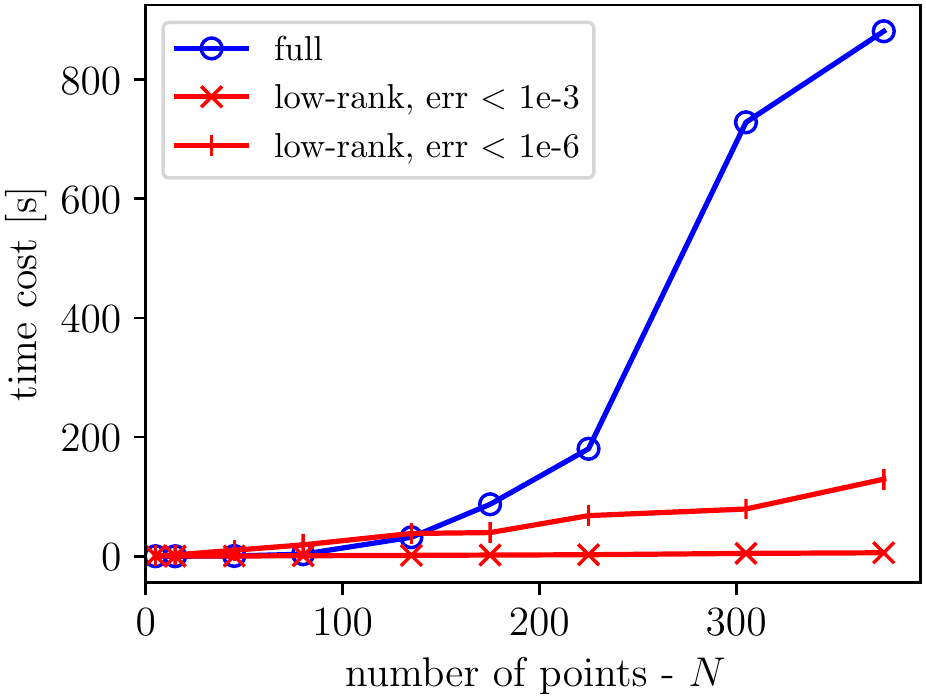}}  \newline
\subfloat[2D, S, anisotropic]{\includegraphics[width=0.49\linewidth]{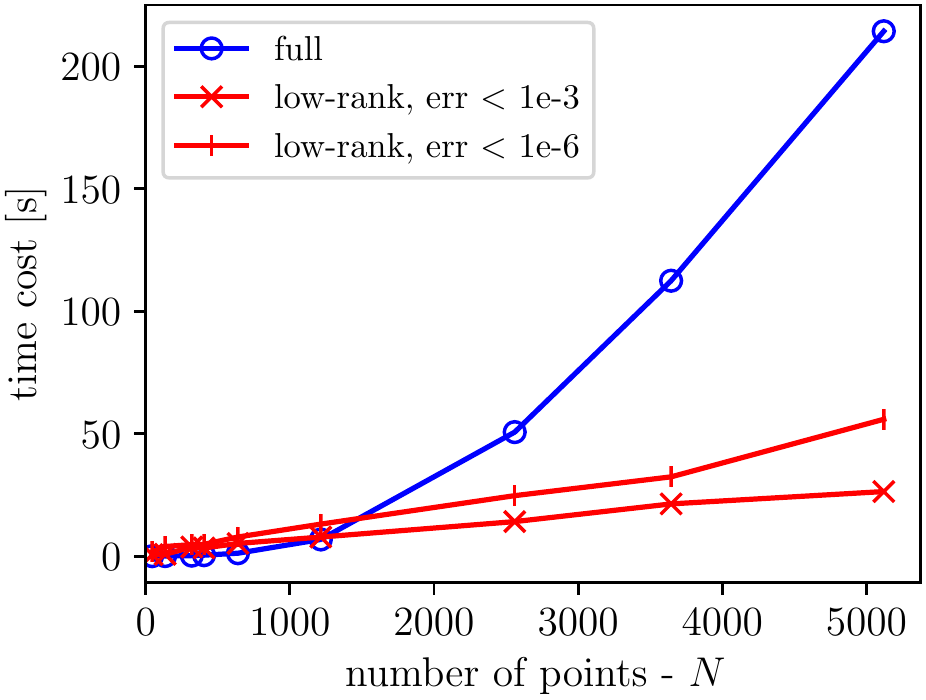}}
\subfloat[3D, S, anisotropic]{\includegraphics[width=0.49\linewidth]{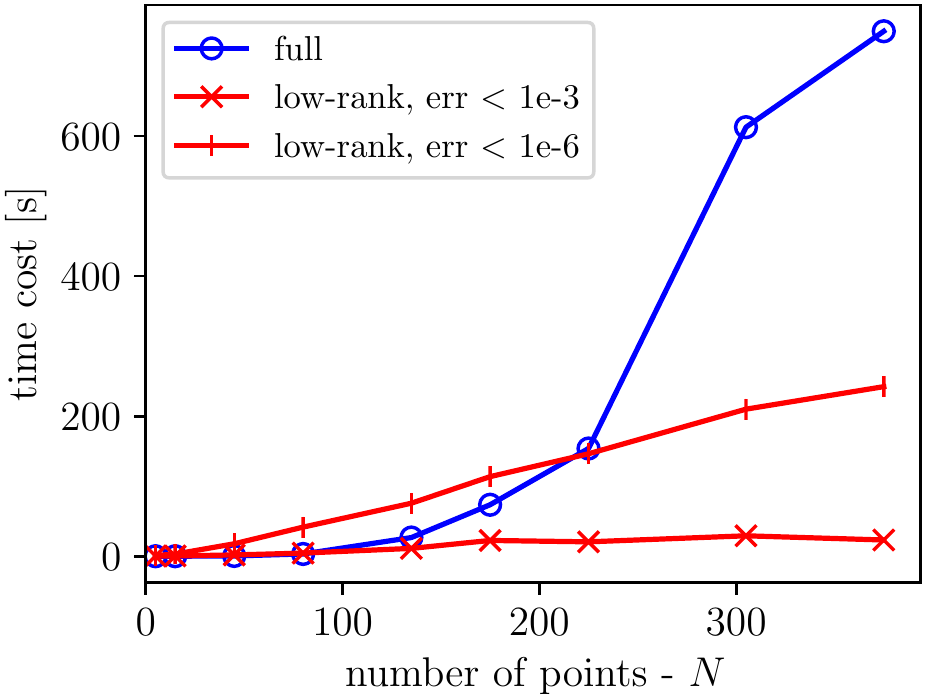}} 
\caption{  The CPU time of GaNi solver  to solve the problem with  isotropic (1st row) and anisotropic (2nd row) material S. Both the full and low-rank solution are computed on the same grid.  The stopping criterion of conjugate gradients for the full solver was set to $10^{-6}$ on the norm of residuum. The minimal residuum iteration for the low-rank approximation was computed with various rank to achieve a required error tolerance defined in \eqref{eq:rel_error}.
}
\label{fig:time_material_2_4}
\end{figure}

\begin{table}[H]
\centering
 \begin{tabular}{l|rrrr|rrrrr}
  \hline \hline
    &\multicolumn{4}{ c|}{ 2D, \squar  } &\multicolumn{5}{ c }{ 3D, \squar }\\ 
  \hline \hline
  $N$ & 45 &  135 &  405 & 1215  & 5 & 15 & 45 & 135 & 175 \\
  $r$ (isotropic) & 3 & 3 & 5  & 7  & 3 & 3 & 3 & 5& 5 \\  
 $r$ (anisotropic) & 5 & 11  & 21 &  31 & 3 & 3 & 5 & 11& 11 \\ 
  \hline \hline
 \end{tabular}
 \caption{Rank $r$ of low-rank solutions that reach the same accuracy as the full solution  for various values of $N$, for isotropic and anisotropic material \squar. The full solver has been computed with grid size $(N,\dotsc,N)$ while sparse solver with $(3N,\dotsc,3N)$.
   The stopping criterion of conjugate gradients for the full solver was set to $10^{-6}$ on the norm of residuum.}
 \label{tab:rank_sparse_solver}
\end{table} 
 

As we are using several low-rank formats and several operations on them, the computational complexities and memory requirements are summarised in Tables~\ref{tab:complexities} and \ref{tab:memory}.
The memory requirements of the FFT-based systems are controlled by memory requirements for material coefficients, preconditioner, solution vector, and possibly other vectors needed to store as a requirement of the linear solver.
Provided that the ranks are kept small, the memory of low-rank solvers scales linearly with $N$, while full solver scales with $N^d$, which makes the method effective particularly for tensor with high order.

  We compare the CPU time of full and low-rank solvers for homogenisation with exact integration (Ga) on the same level of accuracy measured by the energetic norm.
It is achieved by the following procedure.
The reference full solution was computed using conjugate gradient method on the regular grid $(N,\dotsc,N)$ with the tolerance $10^{-6}$ on the norms of residua. 
 In order to achieve the same accuracy as the full solution, the low-rank solver was run on a bigger grid $(\alp N,\dotsc,\alp N)$ with the multiplier $\alp=3$. 
 The rank of low-rank approximations  was increased step-by-step until it achieved a required error tolerance defined in \eqref{eq:rel_error}.
 The iterations of the low-rank solver (for a given rank) are stopped  when the residuum fails to decrease.
 This procedure, which creates a great possibility for a rank reduction in the low-rank solution, is applicable only for problems that allow an exact integration of material coefficients (here material $\squar$). 

The results in Figure~\ref{fig:time_material_0_3} shows that the CPU time scales as $N^d$ for a full solution, and almost linearly for low-rank solutions on the isotropic material $\squar$. In the anisotropic cases the time costs of low-rank solutions are relatively higher but still cheaper than that of the full solution. The difference in these two cases is due to the different ranks of the low-rank solutions. For isotropic material  $\squar$, the solution rank increases only slowly with $N$ , while for its anisotropic counterpart  the rank increases at a faster rate (as tabulated in the Table~\ref{tab:rank_sparse_solver}).
In general, the results show that the low-rank solver are significantly faster for larger $N$, despite being run on a larger computational grid.

We also did the comparison of both solvers for the isotropic and anisotropic material S. However, the smooth material S is better suited for the homogenisation with the numerical integration (GaNi), see the comparison in \cite{VoGe2017FFTHvsFEM}.
Therefore, the comparison of the solvers is run on the same discretisation grid. The ranks of the low-rank solution are chosen such that it achieves a relative error (as defined in \eqref{eq:rel_error}) below $10^{-3}$ or $10^{-6}$.
The ranks remain stable when $N$ increases,  which  makes the CPU time of the low-rank solver almost linear in $N$, as shown in Figure~\ref{fig:time_material_2_4}. 

\section{Conclusion}
This paper is focused on the acceleration of Fourier--Galerkin methods using low-rank tensor  approximations for  spatially  $2$-dimensional and $3$-dimensional problems of numerical homogenisation.
The efficiency of this approach builds on incorporation of the fast Fourier transform (FFT) and low-rank tensor approximation into the iterative linear solvers. 
The computational complexity is reduced to be quasilinear in the size of the discretisation and linear in  spatial dimension  $d$, since on a low-rank tensor of order $d$, the $d$-dimensional FFT can be  performed  as a series of one-dimensional FFTs.
In this paper three formats --- canonical polyadic (CP), Tucker, and tensor train (TT) --- have been considered, and all of them show similar advantage in saving the computational cost.

The main results are summarised as the following:
\begin{itemize}

\item The incorporation of low-rank tensor approximations lead to a significant reduction of memory and computational cost in the solution of the homogenisation problems.

\item The method is more suitable for material coefficients with relatively smaller rank. The low-rank approximation solvers computationally benefits from the better asymptotic behaviour, see Table~\ref{tab:complexities} and \ref{tab:memory}.
The advantage is accentuated for problems of a higher  spatial  dimension $d$ leading to tensors with order $d$.

\item The low-rank approximation can be seen as a model order reduction technique. 

\end{itemize}

Since the low-rank approximation provides a significant memory reduction it allows to compute the solution on a finer grid. Therefore, the proposed method based on low-rank approximation may provide more accurate solution than the conventional method based on full tensors, especially when the material is of a relatively small rank.

\appendix
\section{Low-rank tensor approximations}\label{sec:low-rank-tensor-approximations}
Here we provide more details of the low-rank tensor approximations techniques utilized in this paper. This includes the approximation in CP, Tucker and tensor train formats.

\subsection{The canonical polyadic format}
A canonical polyadic (CP) or $r$-term representation $\MB v_r $ of a tensor $\MB v \in \mathbb K^{N_1 \times \cdots \times N_d}$ ( $\mathbb K$ is either $\mathbb R$ or $\mathbb C$) is a sum of $r$ rank-$1$ tensors, i.e.
\begin{align}
\label{eq:cano}
\MB v \approx \MB v_r  = \sum^r_{i=1} \M{c}[i] \bigotimes^d_{j=1}\bas{j}[i]
\end{align}
with $\bas{j} \in \sK^{r\times N_j}$ and $\bigotimes$ denotes tensor product. 
This format has linear storage size $r \sum^d_{j=1}N_j$. But for $d \geq 3$ and a given $r$, the construction of an error minimizing $\MB v_r$ is not always feasible \cite[Proposition 9.10]{Hackbusch2012book} because the space of CP format tensor with fixed $r$ is not closed \cite[Lemma 9.11]{Hackbusch2012book}. 

\subsubsection{Element-wise multiplication}
The element-wise (Hadamard) product of two tensors of ranks $r$ and $s$ in CP format is computed as:
$$
\MB v_r \odot \MB w_s = \sum^r_{i=1} \sum^s_{k=1} \M{c}_{\MBv}[i]\M{c}_{\MBw}[k] \bigotimes^d_{j=1} \left (\bas{j}_{\MBv}[i] \odot \bas{j}_{\MBw}[k] \right).
$$
This operation has complexity $r s \sum^d_{j=1} N_j $ and the product has a new rank $rs$. 
\subsubsection{Fourier transform}
Due to the linearity and tensor structure of the Fourier transform $\FFT{\VN}$ of a size $\VN\in\sN^d$, a $d$-dimensional Fourier transform of a CP tensor is broken down to a series of $1$-d Fourier transform, i.e., 
$$
\FFT{\VN} (\MB v_r ) = \sum^r_{i=1} \M{c}[i] \bigotimes^d_{j=1} \FFT{N_j}(\bas{j}[i] ).
$$
Hence a FFT on a CP tensor has complexity $ drN\log N $.
\subsubsection{Rank truncation}
Operations (e.g. element-wise multiplication) applied on tensors in CP format usually inflate the representation rank. This calls for a truncation to a prescribed rank or error tolerance.

For $d=2$, this reduction is done by rank truncation based on QR decomposition and singular value decomposition(SVD).
Let the matrices $\B B^{(j)} \in \mathbb{K}^{N_j \times r}$ collect the vectors $\{\bas{j}[i]\}^r_{i=1}$ for the $j$-th dimension, we have their re-orthogonalisations $\B B^{(1)} = \B Q^{(1)} \B R^{(1)}$ and $\B B^{(2)} = \B Q^{(2)} \B R^{(2)}$ by QR decompositions.
A SVD $\B R^{(1)} \B R^{(2)} = \B U^{(1)}\B \Sigma (\B U^{(2)})^\top$ facilitates the truncation. Suppose $\B U^{(1)}_k$, $\B U^{(2)}_k$ and $\B \Sigma_k$ are the truncated ones with rank $k \leq r$, the truncated form of the CP representation \eqref{eq:cano} is 
$$
\MB v_k  = \sum^k_{i=1} \M{c}[i] \hbas{1}[i] \otimes \hbas{2}[i]
$$
where $\hbas{1}[i]$, $\hbas{2}[i]$ are the columns of $\B Q^{(1)} \B U^{(1)}_k$, $\B Q^{(2)} \B U^{(2)}_k$ respectively, and $\M{c}[i]$ are the diagonal entries of $\B \Sigma_k$.

For $d \geq 3$, the $k$-rank form could be obtained by numerical error minimizing procedures \cite{Hackbusch2012book}, e.g. Alternative Least-Squares method. But there is no guarantee that the procedures would converge, and if they would, there is no guarantee that they converge to the global optimum. This is due to the non-closedness of the set of rank-$r$ CP tensors with $d \geq 3$.

\subsection{Tucker format }
\label{sec:tucker}
A Tucker format representation (or tensor subspace representation) of a tensor $\MB v \in \mathbb K^{N_1 \times \cdots \times N_d } \in \C V $ is a linear combination of frames (usually orthogonal bases) of the tensor space $\C V$.
Suppose $ \sV = \bigotimes^d_{j=1} \sV^j$, the subspace $\sV^j$ has basis vectors $\{ \bas{j}[i_j]\in\sK^{N_j}: 1 \leq i_j \leq r_j \}$ with ranks $\Vr=(r_1,\dotsc,r_d)$.
The tensors $\bigotimes^d_{j=1} \bas{j}[i_j]$ for all $1 \leq i_j \leq r_j$ form the bases of the space $\C V$. Then we have a unique coefficient $\M c[i_1, i_2, \dots, i_d] $ for every $\MB v \in \C V$ such that
\begin{align}
\label{eq:tucker}
\MB v \approx \MB v_{\vek r} 
= & \sum^{r_1}_{i_1=1} \cdots \sum^{r_d}_{i_d=1} \MB{c}_{\MB v}[i_1, i_2, \dots, i_d] \bigotimes^d_{j=1} \bas{j}_{\MB v}[i_j],
\end{align}
where $\MB c \in \mathbb K^{r_1 \times \cdots \times r_d }$ is called the core tensor.
Given any prescribed rank vector $\vek r$, an error minimizing approximation $\MB v_{\vek r}$ can be found by a {\em high-order singular value decomposition} (HOSVD) \cite{Lathauwer2000}.
 When the vectors $\{ \bas{j}[i_j]\in\sK^{N_j}: 1 \leq i_j \leq r_j \}$ form only a frame of the subspace $\C V$ (e.g. after addition of two tensors), the core tensor is not unique, however, a representation with orthogonal bases can be obtained by applying QR decomposition to the frames and HOSVD to the accordingly updated core. 

\subsubsection{Element-wise multiplication}
Let another Tucker format tensor with rank $\vek s$ be defined as
\begin{align*}
\MB w_{\vek s} 
 =  \sum^{\vek s}_{\vek k }  \MB c_{\MB w}[\vek k] \bigotimes^d_{j=1}\bas{j}_{\MB w}[k_j]
\end{align*}
the element-wise (Hadamard) product of $\MB v_{\vek r}$ and $\MB w_{\vek s}$ has also a Tucker format
$$
\MB v_{\vek r} \odot \MB w_{\vek s} = \sum^{\vek t}_{\vek l }  \MB c[\vek l]  \bigotimes^d_{j=1} \bas{j}[l_j] 
$$
where $\vek t = \vek r \odot \vek s$ and $\MB c= \MB c_{\MB v} \otimes \MB c_{\MB w}$, i.e. the Kronecker product of the two coefficient tensors. So for any $1 \leq j \leq d$, the index $l_j$ is related to $i_j$ and $k_j$ by 
$
l_j = \overline{i_j k_j} = i_j r_j+k_j 
$,
and $u$ is obtained from $v$ and $w$ through $$ u^{(j)}_{l_j} = u^{(j)}_{\overline{i_j k_j}} = v^{(j)}_{i_j} \odot w^{(j)}_{k_j} \mbox{\;\;\; for \;\; } 1 \leq i_j \leq r_j, \; 1 \leq k_j \leq s_j
$$
Let $N = \max_i N_i$, $r = \max_i r_i$ and $s = \max_i s_i$, 
the computational complexity of the element-wise product is bounded by $d N r s + r^d s^d$, in which the first term is the cost for computing $\{ u^{(j)}_{l_j}: 1 \leq l_j \leq R_j \}_{j=1}^d$, and the second for the Kronecker product of coefficient tensors.


\subsubsection{Fourier transform}
The Fourier transform of $\MB v_{\vek r}$ is
$$
\FFT{\VN} (\MB{v}_{\vek r} ) = \sum^{\vek r}_{\vek i } \MB c [\vek i] \bigotimes^d_{j=1} \FFT{N_j} ( \bas{j}[i_j]) 
$$
which only involves the basis vectors. If FFT is applied, the complexity is of order $\mathcal O(d r N \log N)$.

\subsubsection{Rank truncation}
The Tucker representation \eqref{eq:tucker} can be obtained either by a HOSVD applied on a full tensor or by an operation (e.g. element-wise multiplication) over other Tucker operands. In the first case, an error minimizing rank truncation is readily available due to the property of HOSVD:
$$\sigma^{(j)}_1 \geq \sigma^{(j)}_2 \geq \cdots \geq \sigma^{(j)}_{r_j},\quad\mbox{for } j=1,\cdots,d,$$
where $\sigma^{(j)}_{i_j}$ is the $2$-norm of the $i_j$-th slice of the core tensor $\MB c$ cut on the $j$-th dimension. If the truncation rank is $k_j < r_j$, the error of the truncated representation $\MB v_{\vek k}$ is bounded by
$$
\|\tt v_{\vek r} - \tt v_{\vek k} \| \leq \left [ \sum^d_{j=1} \sum^{r_j}_{i = k_j+1} (\sigma^{(j)}_i)^2 \right ]^{1/2}.
$$

In the second case the bases $\{\bas{j}[i_j] \}^{r_j}_{i_j=1}$ have to be re-orthogonalised first, and then a HOSVD of the updated core tensor is to be made to facilitate the truncation as in the first case. This procedure \cite[as detailed in]{Hackbusch2012book} is analogues to the re-orthogonalisation and SVD for the 2D CP format representations, but  with higher tensor  order.

\subsection{Tensor train format}\label{sec:tensor-train-format}
A tensor train(TT) representation \cite{Oseledets2011} of a tensor $\MB v \in \mathbb K^{N_1 \times \cdots \times N_d}$ can be expressed as a series of consecutive contractions of tensors $\bas{j} \in \mathbb K^{r_{j-1} \times N_j \times r_j}$ of order $3$ for $j=1,\cdots,d$, which are the \emph{carriages} of the tensor train. An equivalent expression in the form of tensor products is

\begin{align*}
\MB v \approx  \MB v_{\vek r} &= 
\sum^{r_1}_{i_1=1}\cdots \sum^{r_{d-1}}_{i_{d-1}=1} \bas{1}_{\MB v}[1,:,i_1] \otimes \bas{2}_{\MB v}[i_1,:,i_2]  \otimes \cdots \otimes \bas{d}_{\MB v}[i_{d-1},:,1]
\end{align*}
$\vek r$ is the TT-rank of $\MB v$ with a constrain $r_0 = r_d =1$ to keep the elements of $\MB v$ scalars. The TT format is stable in the sense that for any prescribed $\vek r$ an error minimizing $\MB v_{\vek r}$ can always be constructed by a series of SVDs on consecutive matricisations of $\MB v$.

\subsubsection{Element-wise multiplication}
Let another TT format tensor with rank $\vek s$ be defined as
$$
\MB w_{\vek s} = 
\sum^{s_1}_{i_1=1}\cdots \sum^{s_{d-1}}_{i_{d-1}=1} \bas{1}_{\MB w}[1,:,i_1] \otimes \bas{2}_{\MB w}[i_1,:,i_2]  \otimes \cdots \otimes \bas{d}_{\MB w}[i_{d-1},:,1]
$$
with $\bas{j}_{\MB w} \in \mathbb K^{s_{j-1} \times N_j \times s_j}$. The element-wise  product of $\MB v_{\vek r}$ and $\MB w_{\vek s}$  can also be expressed in TT format:
$$
\MB v_{\vek r} \odot \MB w_{\vek s} =
\sum^{t_1}_{i_1=1}\cdots \sum^{t_{d-1}}_{i_{d-1}=1} \bas{1}[1,:,i_1] \otimes \bas{2}[i_1,:,i_2]  \otimes \cdots \otimes \bas{d}[i_{d-1},:,1]
$$
where $\vek t = \vek r \odot \vek s$ and $\bas{j} = \bas{j}_{\MB v} \ast \bas{j}_{\MB w}$. 
Here the $\ast$ denotes one type of Khatri–Rao product \cite{Khatri1968} which makes Kronecker product only in the first and third dimensions, i.e. it yields an order $3$ tensor  $\bas{j} \in \mathbb K^{r_{j-1} s_{j-1}\times N_j \times r_{j} s_{j}}$. The complexity of the element-wise product is of order $\mathcal O(d N r^2 s^2)$ with $N$, $r$ and  $s$ as defined in the  subsection~\ref{sec:tucker}.

\subsubsection{Fourier transform}
The Fourier transform of $\MB v_{\vek r}$ can also be carried out by doing $1$-D transforms on each \emph{carriage}: 
$$
\FFT{\VN} (\MB v_{\vek r} ) =
\sum^{r_1}_{i_1=1}\cdots \sum^{r_{d-1}}_{i_{d-1}=1} \FFT{N_1} (\bas{1}_{\MB v}[1,:,i_1])
\otimes 
\FFT{N_2} (\bas{2}_{\MB v}[i_1,:,i_2])  \otimes \cdots \otimes 
\FFT{N_d} (\bas{d}_{\MB v}[i_{d-1},:,1])
$$
in which the $\FFT{N_j}(\cdot)$ is made on the fibres along the second mode. If FFT is applied here, the number of operations is of order $\mathcal O(d r^2 N \log N)$.

\subsubsection{Rank truncation}
The tensor train representation \eqref{eq:TT} can be obtained either by transforming a full tensor into tensor train format by using $d-1$ sequential SVDs applied on auxiliary matrices of the tensor (known as TT-SVD) \cite{Oseledets2011}, or as a result of operations (e.g. additions or multiplications) over tensor train operands.
In the first case, an error minimising rank truncation could be directly carried out in the TT-SVD process.
The truncation has an error bound $(\sum^{d-1}_{k=1} \epsilon^2_k)^{1/2} $, where $\epsilon_k$ is the Frobenious norm error introduced by the truncation of the $k$-th SVD.
In the second case, a re-orthogonalisation has to be done in the first place, this is followed by $d-1$ sequential SVDs on unfolded \emph{carriages}. This process is known as TT-truncation (also called rounding).

For the first case, the complexity of truncation is the same as that for the TT-SVD, which is of order $\C O(N^{d+1})$. A cheaper alternative for TT-SVD is TT-cross approximation as introduced in \cite{Oseledets2010}. The complexity of TT-truncation in the second case is of order $\C O(dNr^3)$.

\subsection*{Acknowledgements}
Funded by the Deutsche Forschungsgemeinschaft (DFG, German Research Foundation) --- project number MA2236/27-1.

 Martin Ladeck\'{y} was supported by the Czech Science Foundation through projects No. GA\v{C}R 17-04150J, by the
 Center of Advanced Applied Sciences (CAAS),
 financially supported by the European Regional Development Fund (through project No.\ CZ.02.1.01/0.0/0.0/16\_019/0000778), and by the Grant Agency of the CTU in Prague, grant No. SGS19/002/OHK1/1T/11.


 

\end{document}